\documentclass{aims}
\usepackage{amsmath}
  \usepackage{paralist}
  \usepackage{caption, subcaption}

  \usepackage{graphics} 
  \usepackage{epsfig} 
\usepackage{graphicx}  \usepackage{epstopdf}
 \usepackage[colorlinks=true]{hyperref}
   \definecolor{OliveGreen}{rgb}{0,0.6,0}
\hypersetup{urlcolor=OliveGreen, citecolor=OliveGreen, linkcolor=OliveGreen}

\def\currentvolume{X}
 \def\currentissue{X}
  \def\currentyear{200X}
   \def\currentmonth{XX}
    \def\ppages{X--XX}
     \def\DOI{10.3934/xx.xx.xx.xx}

\newtheorem{theorem}{Theorem}[section]

\newtheorem{lemma}[theorem]{Lemma}

\theoremstyle{definition}
\newtheorem{definition}[theorem]{Definition}

\usepackage{amsmath,amssymb}
\usepackage{empheq}
\usepackage{nicefrac}
\usepackage{bbm}

\title[Macroscopic descriptions of follower-leader systems] 
      {Macroscopic descriptions of follower-leader systems}

\author[S. Bernardi, G. Estrada-Rodriguez, H. Gimperlein and K. J. Painter]{}


 \email{sara.bernardi@polito.it}
 \email{estradarodriguez@ljll.math.upmc.fr}
 \email{h.gimperlein@hw.ac.uk}
 \email{kevin.painter@polito.it}

\thanks{G. E. R. was supported by The Maxwell Institute Graduate School in Analysis and its Applications, a Centre for Doctoral Training funded by the UK Engineering and Physical Sciences Research Council (grant EP/L016508/01), the Scottish Funding Council, Heriot-Watt University and the University of Edinburgh.}

\thanks{$^*$ Corresponding author: Heiko Gimperlein}

\begin{document}
\maketitle

\centerline{\scshape Sara Bernardi}
\medskip
{\footnotesize
  \centerline{Department of Mathematical Sciences ``G. L. Lagrange'', Politecnico di Torino }
   \centerline{Corso Duca degli Abruzzi 24, 10129 Torino, Italy}}

\medskip

\centerline{\scshape Gissell Estrada-Rodriguez}
{\footnotesize
 \centerline{Laboratoire Jacques-Louis Lions, Sorbonne-Universit\'{e}, 4, pl. Jussieu, F-75005 Paris, France}
} 

\medskip

\centerline{\scshape Heiko Gimperlein$^{*}$}
\medskip
{\footnotesize
 \centerline{Maxwell Institute for Mathematical Sciences and Department of Mathematics, Heriot-–Watt University}
   \centerline{Edinburgh, EH14 4AS, United Kingdom}
   \centerline{and Institute for Mathematics, University of Paderborn, Warburger Str.~100, 33098 Paderborn, Germany }
} 

\medskip

\centerline{\scshape Kevin J. Painter}
\medskip
{\footnotesize
 \centerline{Interuniversity Department of Regional and Urban Studies and Planning, Politecnico di Torino}
   \centerline{Torino, 10125, Italy}
}

\bigskip

 \centerline{(Communicated by the associate editor name)}

\begin{abstract}
The fundamental derivation of macroscopic model equations to describe swarms based on 
microscopic movement laws and mathematical analyses into their self-organisation 
capabilities remains a challenge from the perspective of both modelling and analysis. In this paper we 
clarify relevant continuous macroscopic model equations that describe follower-leader 
interactions for a swarm where these two populations are fixed. We study the behaviour 
of the swarm over long and short time scales to shed light on the number of leaders needed 
to initiate swarm movement, according to the homogeneous or inhomogeneous nature of the 
interaction (alignment) kernel. The results indicate the crucial role played by the 
interaction kernel to model transient behaviour.
\end{abstract}

\section{Introduction}

Collective movements describe the tendency of a group of individuals to coordinate 
their motion in a manner that generates net flow of the entire population. Examples 
range from cells to animals, from migrating cell clusters during development and 
cancer invasion \cite{haeger2015,mayor2016} to shifting bird flocks and fish shoals; the 
latter extend to the kilometre-spanning shoals formed from hundreds of millions of 
herrings \cite{makris2009}. A point of recent interest concerns the potential 
division of a population into ``leaders'' and ``followers'' and, consequently, how 
leaders influence swarm dynamics. Clearcut leadership could result from 
experience, age or prior knowledge: the presence of older birds 
improves path efficiency in migrating cranes \cite{mueller2013}; post-menopausal orcas adopt leading 
positions during pod foraging \cite{brent2015}; only the ``househunting'' scouts 
know the final destination of a new honeybee nest \cite{seeley}. Subtler leaders
arise within superficially identical populations, such as the presence of 
faster or ``braver'' individuals in fish shoals and bird flocks \cite{reebs2000,nagy2010}. 
Leader/follower statuses also occur in a host of cellular systems, ranging from
collective movements of aggregated amoebae to sprouting blood capillaries during 
development, physiology and disease \cite{mayor2016}.

Most theoretical descriptions of collective migration have employed 
agent/particle-based approaches, e.g. \cite{Couzin,Couzin2,Cucker2,Cucker,vicsek}, 
see also \cite{berdahl2018}. Despite the plethora of models, they typically share a  
so-called set of ``first principles of swarming'' \cite{Carrillo}: specifically, 
particle trajectories governed by a combination of repulsion (preventing collisions), 
attraction/cohesion (preventing dispersal) and alignment of direction/velocity
according to neighbour positions, each operating over specified interaction ranges. Models 
based on these features reproduce a wide variety of collective migration 
phenomenologies, e.g. see \cite{berdahl2018}. In the context of follower-leader systems, 
a key finding has been that swarms can be efficiently guided by a small number of 
anonymous (i.e. not clearly distinct from the crowd) informed individuals \cite{Couzin}. 
These individuals influence their neighbours, which in turn influence further followers
and knowledge is relayed through the swarm. Surprisingly, as the swarm population 
increases a diminishing fraction of leaders is needed to achieve the 
same level of guidance efficiency \cite{Couzin}. Follower-leader alignment strategies have also
been incorporated within other swarming studies. For example, in 
\cite{jadbabaie2003coordination} a ``transient leadership'' model was considered to 
imitate bird flocks while ``hierarchical leadership'' was studied in \cite{shen2007cucker} 
within the framework of a Cucker-Smale model.

Beyond agent-based models, a plethora of continuous models have been proposed, 
for example see \cite{bernardi2021,eftimie2012,Eftimie,mogilner1999,painter2015,topaz2006}. In 
common with their individual-based counterparts, movement is governed by a set or 
subset of attracting/repelling/aligning interactions, typically generating
integro-differential equations of parabolic (e.g. \cite{mogilner1999,painter2015}) or hyperbolic (e.g.
\cite{bernardi2021,eftimie2012}) form.  While these models gain analytical tractability, their 
connection to individual-level behaviour is, inevitably, blurred.

The principal objective of the current paper is to clarify a relevant form of partial differential 
equation system to describe follower-leader interactions within swarming populations. Specifically, we focus on swarms that are clearly divided into two
distinct subpopulations: a leader population with knowledge of the destination,
and an unaware follower population. A paradigm for this form of follower-leader
system is provided by honey bee swarms. Prior to swarming, a small proportion of scouts search the environment for a new nest,
eventually arriving at a consensus for the nest location \cite{seeley}. Successful translocation of
the swarm subsequently relies on this small population of knowledgeable bees leading
the uninformed colony to this new location. Scout bees appear to transmit their
information by performing a sequence of nest-directed streaks through the upper
swarm \cite{beekman2006does,seeley,greggers2013}, a behaviour that is believed to accentuate their visibility
to followers. Inevitably, such behaviour quickly brings scouts to the leading edge, so
a shift in behaviour would be necessary to prevent losing contact with the swarm.
One proposal, suggested in \cite{seeley}, is slowly travelling backwards, perhaps along the
bottom or sides to minimise their influence, before streaking again. This suggests
that leaders switch between active and passive states, e.g. as proposed in \cite{bernardi2018particle}, where in the former they attempt to maximise their influence on the swarm orientation and
in the latter it is minimised.

Taking inspiration from this and other follower-leader systems,
we begin with an underlying particle description for movement that is based on a
relatively minimal set of assumptions. Specifically. we propose velocity alignment
interactions that could lead to coherent and guided swarm movement. We subsequently
examine possible hyperbolic and parabolic scaling limits. We discuss the
plausibility of these macroscopic models, their relevance in the context of follower-leader
systems such as the bee swarms described above, and the new modelling
challenges raised. We further present some preliminary simulations, indicating the extent to which the model can describe leader-follower behaviour.

\section{Swarm description}\label{sec: swarm description}

We assume a swarm attempting to reach some target
destination and divided into two main subpopulations: followers ($f$) and leaders.
Leaders are assumed to have knowledge of the target while the followers are completely
uninformed. The leader population is further subdivided into those engaging
as active leaders ($a$) and those behaving passively ($p$). Active leaders are defined as
that fraction of the leader population who are moving in the direction of the target,
and can potentially engage in behaviour designed to communicate directional
information regarding the target. Passive leaders are assumed to move away from
the target, e.g. towards the back of the swarm where they can reengage in active
leader behaviour. We remark that the total number of leaders is considered to be
small, hence we can describe the leaders as a discrete system (it is noted that for bee
populations around 3\%-5\% of the total swarm is believed to be a leader \cite{seeley}, while
for fish a minority of informed individuals is known to lead a shoal to food \cite{reebs2000}).
Transitions are only allowed to occur between passive and active leaders, i.e. there
is no transition between follower and leader type. Transitions will be motivated by
the objective of maintaining contact with the swarm, i.e. an active leader that has
moved to the leading of edge of the swarm changes behaviour to avoid outruning
the swarm, while a passive leader moving backwards similarly switches back. Notably, though, we do not explictly incorporate an ``attractive'' or ``cohesive'' response into individual behaviours, i.e. where individuals are attracted towards higher swarm density regions.

An obviously important quantity in the model is the position of the swarm with respect to the destination or target. Given a generic member located at $\mathbf{x}_k$, its distance to the target 
at $\mathbf{x}_{\textup{target}}$ is denoted by $I_{\textnormal{target}}(\mathbf{x}_k(t))=|\mathbf{x}_k(t) - \mathbf{x}_{\textup{target}}|$. We fix a corresponding field $\mathbf{b} \simeq \nabla I_{\textnormal{target}}$ within the domain outside the target, such that $- \mathbf{b}$  defines the direction for the leaders to the target. Specifically, we choose $\mathbf{b}=\nabla I_{\textnormal{target}}$ to be divergence free at the target and tangent to the boundary of the domain. As $|\mathbf{x}| \to \infty$, we impose $\mathbf{b}(\mathbf{x}) \to \nabla I_{\textnormal{target}}$.

Following the approach in \cite{alt1980biased}, we define $\bar{\sigma}_i (\mathbf{x},t,\theta,\tau)$, for $i \in \left\{f,p,a\right\}$, as the {\em microscopic densities} of followers, passive leaders and active leaders at position $\mathbf{x}$, time $t$, moving in direction $\theta$ for some time $\tau$. Integrating with respect to the microscopic quantity $\tau$, we 
denote $\sigma_i (\mathbf{x},t,\theta) = \int_0^t \bar{\sigma}_i (\mathbf{x},t,\theta,\tau)d\tau$
as the density of each subpopulation at position $\mathbf{x}$, time $t$ and moving in 
direction $\theta$ and we call it a \emph{mesoscopic density}. Subsequently integrating over $\theta$ generates the {\em macroscopic densities}
\begin{equation}
\rho_i (\mathbf{x},t) = \int_0^t \int_S \bar{\sigma}_i(\cdot,\theta,\tau)d\theta d\tau\,,\label{eq: macroscopic definition}
\end{equation}
where $S=\{\mathbf{\theta}\in \mathbb{R}^n:\ |\mathbf{\theta}|=1\}$ is the unit sphere in $\mathbb{R}^n$ {and the measure $d\theta$ is normalized so that the surface area $|S|=1$}. Dropped subscripts are used to denote total population densities, i.e.
$\bar{\sigma} =\bar{\sigma}_f+\bar{\sigma}_p+\bar{\sigma}_a$, $\sigma =\sigma_f+\sigma_p+\sigma_a$ 
and $\rho =\rho_f+\rho_p+\rho_a$. 

Each subpopulation follows a random walk process, described according to a set of rules minimally chosen to describe a generic swarm behaviour but restricting excessive complexity. Motivation behind some of these choices can be found by paradigm examples, such as bee swarms, although we are not aiming to describe a
specific system.

Follower behaviour is described by the following set of assumptions.
    \begin{enumerate}
        \item[1F.] Trajectories comprise of straight line motions interrupted by (effectively) instantaneous reorientations, where the new direction of motion is randomly chosen. This movement is called a velocity-jump process \cite{othmer1988}. Individuals stop (i.e. reorient) with a rate given by a fixed parameter $\beta$.
        \item[2F.] At each reorientation, with probability $\zeta\in(0,1)$ it selects a new direction of motion $\eta$, taken to be symmetrically distributed with respect to the previous one according to 
        \[k(\mathbf{\theta};\mathbf{\eta})=\tilde{k}(|\eta-\theta|)\ .\]
        Because $\tilde{k}$ is a probability distribution, it is normalized to $\int_S\tilde{k}(|\theta-e_1|)d\theta=1$ where $e_1 = (1, 0, . . . , 0)$. The turn angle operator $T: L^2(S) \rightarrow L^2(S)$ is defined as
        \begin{equation}
    T\phi(\eta)=\int_Sk(\theta;\eta)\phi(\theta)d\theta\ .\label{eq: turning angle simple}
\end{equation}
Trivially, $\tilde{k}$ could be a uniform distribution; more generally, a bias according to the previous orientation would incorporate an element of persistence of orientation, so we consider $\tilde{k}(|\eta-\theta|)$.
        \item[3F.] With probability $1-\zeta$ the follower instead aligns with the orientation of the local population, which is defined by
        \begin{equation}
        \Lambda(\mathbf{x},\theta,t)=\frac{\mathcal{J}(\mathbf{x},t)}{|\mathcal{J}(\mathbf{x},t)|}\ ,\ \mathcal{J}(\mathbf{x},t)=\int_{\mathbb{R}^n}\int_SK(|\mathbf{y}-\mathbf{x}|)\sigma(\mathbf{y},t,\theta)\theta d\mathbf{y}d\theta\ .\label{eq: mean direction}
        \end{equation}
        Here $K$ is an interaction kernel. { If the flux $\mathcal{J}(\mathbf{x},t)=0$, we assume that $\Lambda(\mathbf{x},\theta,t)$ takes the value $\theta$ \cite{degond2008continuum}.} 
        A generalisation of the above would be to choose something of the form:
        \begin{equation}
        \Lambda(\mathbf{x},\theta,t)=\frac{\mathcal{J}(\mathbf{x},t)}{|\mathcal{J}(\mathbf{x},t)|}\ ,\ \mathcal{J}(\mathbf{x},t)=\int_{\mathbb{R}^n}\int_SK(|\mathbf{y}-\mathbf{x}|)\left(\lambda_f\sigma_f+\lambda_p\sigma_p+\lambda_a \sigma_a\right) d\mathbf{y}d\theta\ ,\label{eq: mean direction alt form}
        \end{equation}        
        where $\lambda_{f,p,a} \ge 0$ would reflect the capacity
of followers to differentiate between the various subpopulations during its
orientation. A choice $\lambda_a \gg \lambda_{f,p}$ would assume active leaders are considerably
more conspicuous than the other populations and dominate the orientation
behaviour, for example through engaging in specific movement behaviour
or issuing vocal commands; a choice $\lambda_p$ = 0, on the other hand, would assume passive leaders make themselves ``invisible'' to the followers. We note that the limitations of the above homogeneous choice will be discussed in Section \ref{sec:macrodiscussion}. In particular, an inhomogeneous variant, given by $\Lambda(\mathbf{x},t)=\nu\mathcal{J}(\mathbf{x},t)$, where $\nu$ is a relaxation frequency that depends on the norm of $\mathcal{J}(\mathbf{x},t)$, will subsequently be considered in Section \ref{sec: inhom}.
        \item[4F.] Followers move with a fixed speed $c_f$.
    \end{enumerate} 

\noindent

For the leaders dynamics, we simplify the bias random walk process by assuming the turning distribution to be highly concentrated into choosing a specific movement: (i) a flight towards the target for the active leaders, (ii) a flight in the opposite direction for passive leaders. The limiting case of small randomness follows the discrete model proposed in \cite{bernardi2018particle}, and is predicated on the assumption that leaders have enough knowledge of the target to produce highly directed movements

\begin{figure}
    \centering
    \includegraphics[scale=0.4]{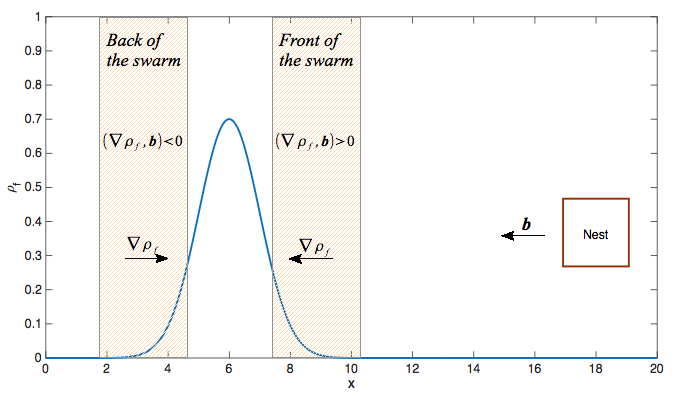}
    \caption{Illustration of switching between streakers and passive leaders.}
    \label{fig: swarm1}
\end{figure}

{
The dynamics of passive leaders are then specified by the following rules: }
    \begin{enumerate}
        \item[1P.] No alignment according to other individuals.
        \item[2P.] Passive leaders move in the direction $\mathbf{c}$, in the opposite direction to the target.
        
        \item[3P.] Passive leaders move with a speed $c_p$.
        \item[4P.] Transitions between active and passive leader take place at the front and rear edges of the swarm, as described below.
    \end{enumerate}
Finally, active leaders move according to:
\begin{enumerate}
   \item[1A.] Given their knowledge of the target, they move as ballistic particles in direction $-\mathbf{b}$ (towards the front of the swarm as represented in Figure \ref{fig: swarm1}), with a fixed maximum speed given by $c_a$. Thus, no randomness is assumed.
   \item[2A.] Transitions between active and passive leader take place at the front and rear edges of the swarm, as described below.    
   \item[3A.] Active leaders move with a speed $c_a$.
\end{enumerate}
We note that movements into orientations outside $-\mathbf{b}$ are not permitted for active leaders within this model. Microscopic densities can then be inferred as singular distributions, where the local macroscopic density is concentrated into orientation $-\mathbf{b}$: i.e., 
$\sigma_a(\mathbf{x},t,\theta) = \rho_a(\mathbf{x},t) \delta (\theta+\mathbf{b})$ where $\delta$ is the Dirac delta function. 

\begin{figure}
    \centering
    \includegraphics[scale=0.25]{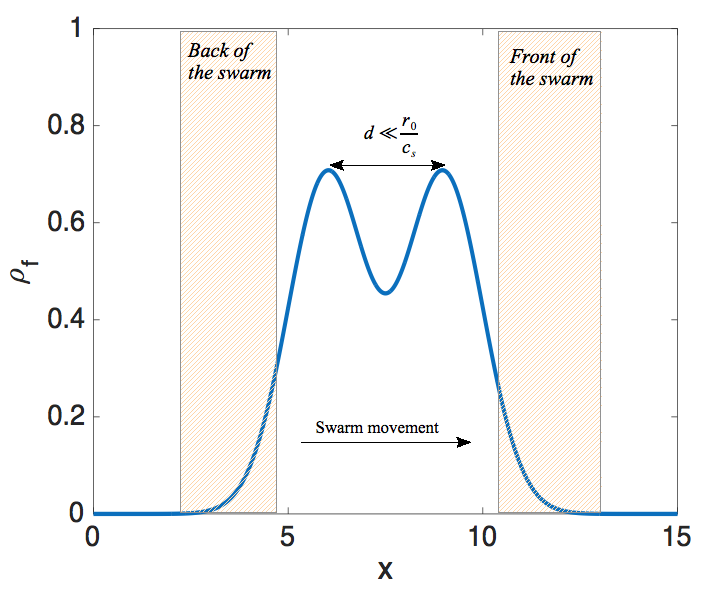}\hspace{0.2cm}\includegraphics[scale=0.25]{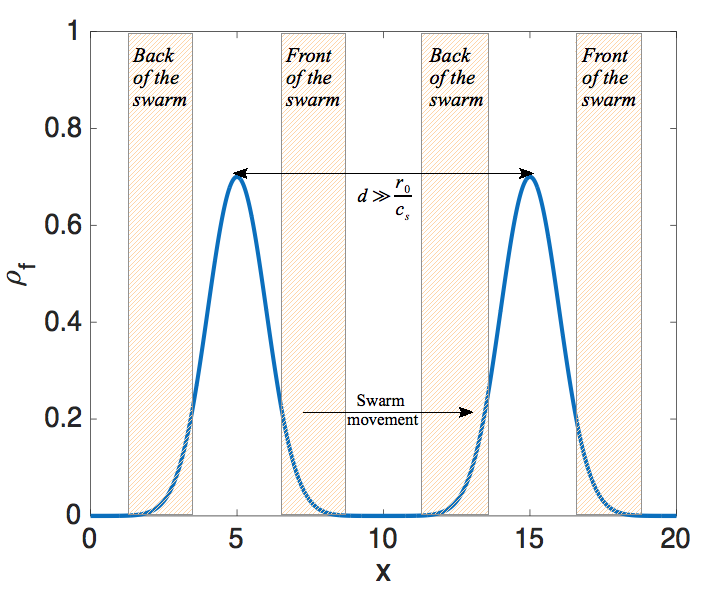}
    \caption{Illustration of different swarm shapes.}
    \label{fig: swarm23}
\end{figure}

Suppose a Gaussian-type curve for the swarm profile. Physically, the transition rate from active to passive leader, $R_{ap}$, should be localised to the front edge of a swarm: in a region where 
$(\nabla\rho_f,\mathbf{b})>0$, as illustrated in Figure \ref{fig: swarm1}. As follower densities 
decrease, i.e. $\rho_f\rightarrow 0$, $R_{ap}$ should be bounded away from zero: $R_{ap}\geq r_0>0$ as 
$\rho_f\rightarrow 0$, where $r_0$ is the minimal conversion rate outside the swarm. 
The rules for transition rate from passive to active leaders, $R_{pa}$, follow a similar description.
Transitions should be concentrated to the rear edge of the swarm, i.e. where $(\nabla\rho_f,\mathbf{b}) <0$. Again, as follower densities decrease we require $R_{pa}\geq r_0>0$.

For more complicated swarm profiles, e.g. dual peaked as in Figure \ref{fig: swarm23}, let $c_s/r_0$ define a length scale on which the density of active leaders (and hence also that of passive leaders) 
varies. If the spacing between the two peaks is sufficiently larger than $c_s/r_0$ then active leaders
consider the blobs as separate swarms and there are separate fronts and backs for each swarm. For 
peaks spaced below this distance, the majority of active leaders fly towards the leading edge 
before converting to passive leaders and it can be considered a single swarm with a single leader population.

We note that while a change from passive to active leader involves transition into a population with fixed orientation, an active to passive change demands transition into a population whose orientation is distributed over the unit sphere. We define $R_{ap}=\int_S\bar{R}_{ap}(\theta)d\theta$ and similarly for $R_{pa}$.

\subsection{Microscopic description}

Following the swarm description given in Section \ref{sec: swarm description}, the 
population densities satisfy systems of integro-differential equations as follows. For the active leaders we consider 
\begin{equation}
    \begin{aligned}
    \begin{cases}
    & \left(\partial_\tau+\partial_t-c_a\mathbf{b}\cdot\nabla \right)\bar{\sigma}_a(\cdot,\tau)=-R_{ap}\bar{\sigma}_a(\cdot,\tau)\ ,\\ &\bar{\sigma}_a(\cdot,\tau=0)=R_{pa}\int_0^t\int_S\bar{\sigma}_p(\cdot,\theta,\tau)d\theta d\tau\ ,\label{eq: active leaders}
    \end{cases}
    \end{aligned}
    \end{equation}
    where $(\cdot)$ denotes the space and time dependence $(\mathbf{x},t)$. The left hand side of the first equation describes the movement of the individuals towards the nest, for time $\tau$, while the right hand side gives the density of active leaders that change to passive. The density of active leaders that start a new run, at $\tau=0$, is given by passive leaders switching to active. This is described by the second equation.
    
    {Note that the vector field $(\partial_t, \partial_\tau, c_a\mathbf{b}\cdot\nabla)$ in \eqref{eq: active leaders} is transversal to the hypersurface $\{\tau = 0\}$, so that the second equation provides the necessary initial conditions.}  
    
The passive leaders obey the following system
\begin{equation}
    \begin{aligned}
    \begin{cases}
    &\left(\partial_{\tau}+\partial_t+c_p\mathbf{c}\cdot\nabla \right)\bar{\sigma}_p(\cdot,\tau)=-R_{pa}\bar{\sigma}_p(\cdot,\tau)\ ,\\ &\bar{\sigma}_p(\cdot,\tau=0)={R}_{ap}\int_0^t\int_S\bar{\sigma}_a(\cdot,\theta, \tau)d\theta d\tau\ .\label{eq:governing_eq2new}
    \end{cases}
    \end{aligned}
    \tag{6.a}
    \end{equation}
   Including some randomness in the movement of the passive leaders we can write the system \eqref{eq:governing_eq2new} as follows
\begin{equation}
    \begin{aligned}
    \begin{cases}
    &\left(\partial_{\tau}+\partial_t+c_p\theta\cdot\nabla \right)\bar{\sigma}_p(\cdot,\theta,\tau)=-R_{pa}\bar{\sigma}_p(\cdot,\tau)-\beta\bar{\sigma}_p(\cdot,\theta,\tau)\ ,\\ &{ \bar{\sigma}_p(\cdot,\theta,\tau=0)=\bar{R}_{ap}(\theta)\int_0^t\int_S\bar{\sigma}_a(\cdot,\eta, \tau)d\eta d\tau+B(\theta)\int_0^t\int_S\beta\bar{\sigma}_p(\cdot,\eta,\tau)d\eta d\tau}\ .\label{eq: passive 2}
    \end{cases}
    \end{aligned}
    \tag{6.b}
    \end{equation}
   The passive leaders in this case also move according to a velocity jump process as described in 1F, with the same stopping rate $\beta$. The new direction of motion $\eta$ is given by the turn angle operator {$B(\eta)$. To be more specific, we define $B(\eta)=pT(\eta)+(1-p)T_{\rho_f}$ where $p\in(0,1)$ { and we fix $T_{\rho_f}=(\mathbf{b}\cdot\nabla\rho_f)\nabla\rho_f$}. The first term in $B(\eta)$ describes the random movement where a new direction of motion is taken arbitrarily, in this case it only depends on the final direction. The second term describes a vector field pointing in the direction of the back of the swarm, where the switch between passive and active leaders occurs. The analysis of equation \eqref{eq: passive 2} is postponed to Section \ref{sec: ranodmness passive}.}

For the followers we consider the  approach in \cite{estrada2019interacting}  and then we have
\stepcounter{equation}
\begin{equation}
    \begin{aligned}
    \begin{cases}
    &\left(\partial_{\tau}+\partial_t+c_f\theta\cdot\nabla \right)\bar{\sigma}_f(\cdot,\theta,\tau)=-\beta\bar{\sigma}_f(\cdot,\theta,\tau) \ ,\\
    &\bar{\sigma}_f(\cdot,\eta,\tau=0)  =\int_SQ(\eta,\theta)\int_0^t\beta\bar{\sigma}_f(\cdot,\theta,\tau)d\tau d\theta\ ,\label{eq: alignment first first}
    \end{cases}
    \end{aligned}
\end{equation}
where 
\begin{equation}
Q(\eta,\theta)=\zeta \tilde{k}(|\eta-\theta|)+(1-\zeta)\Phi(\Lambda\cdot\eta)\ .\label{eq: alignment}
\end{equation}
$\Phi(\Lambda\cdot\eta)$ is the distribution of the new aligned direction and satisfies $\int_S\Phi(\Lambda\cdot\eta)d\eta=1$.
Integrating with respect to $\tau$  the systems (\ref{eq: active leaders}) and (\ref{eq:governing_eq2new})
we obtain, for $\rho_a(\cdot)=\int_0^t\bar{\sigma}_a(\cdot,\tau)d\tau$ and $\rho_p(\cdot)=\int_0^t\bar{\sigma}_p(\cdot,\tau)d\tau$, respectively,
\begin{align}
\partial_t\rho_a-c_a\mathbf{b}\cdot\nabla\rho_a &=-R_{ap}\rho_a+R_{pa}\rho_p\ ,\label{eq: first s}\\
    \partial_t\rho_p+c_p\mathbf{c}\cdot\nabla\rho_p  &=-R_{pa}\rho_p+{R}_{ap}\rho_a\ .
\end{align}

Similarly, integrating system (\ref{eq: alignment first first}) with respect to $\tau$ and using the definition of $T$ given in (\ref{eq: turning angle simple}) we obtain 
\begin{equation}
    \partial_t\sigma_f+c_f\theta\cdot\nabla\sigma_f=-\beta\sigma_f+\zeta\beta T\sigma_f+(1-\zeta)\beta\Phi(\Lambda\cdot\theta)\rho_f\ .\label{eq: first f}
\end{equation}

\section{Macroscopic PDE description}\label{sec: macroscopic PDE}

As shown in the previous section, the densities $\rho_a(\mathbf{x},t)$, $\rho_p(\mathbf{x},t)$ and  $\sigma_f(\mathbf{x},t,\theta)$ satisfy the following system of kinetic equations,
\begin{align}
\partial_t\rho_a-c_a\mathbf{b}\cdot\nabla\rho_a & =-R_{ap}\rho_a+R_{pa}\rho_p\label{eq: streakers}\ ,\\ \partial_t\rho_p+c_p\mathbf{c}\cdot\nabla\rho_p  &=-R_{pa}\rho_p+{R}_{ap}\rho_a \label{eq: passive}\ ,\\ \partial_t\sigma_f+c_f\theta\cdot\nabla\sigma_f & = -\beta\sigma_f+\zeta\beta T\sigma_f+(1-\zeta)\beta\Phi(\Lambda\cdot\theta)\rho_f \label{eq: followers}   \ .
\end{align}

We control that the population of the followers (\ref{eq: followers}) is conserved:
\begin{align}
    \partial_t\int_S\sigma_fd\theta+c_f\int_S\theta\cdot\nabla\sigma_fd\theta=-\beta\int_S\sigma_fd\theta+\zeta\beta\int_ST\sigma_fd\theta+(1-\zeta)\beta\rho_f=0\ , \label{eq: conservation followers}
\end{align}
since $\int_S\Phi(\Lambda\cdot\theta)d\theta=1$ and $\int_S Td\theta=1$ by definition.\\

In the following subsections we study the macroscopic behaviour of the above systems by considering different scaling limits. In
Subsection \ref{sub: hyperbolic limit} we study a hyperbolic limit of the system \eqref{eq: active leaders}-\eqref{eq:governing_eq2new}-\eqref{eq: alignment first first}, where the leaders are considered to be discrete. Then, in Subsection \ref{sec: passive leaders with noise}, a hyperbolic limit of system \eqref{eq: active leaders}-\eqref{eq: passive 2}-\eqref{eq: alignment first first}, where a random component of the movement is included for the case of passive leaders. Later, in Subsection \ref{sub: diffusion limit} we consider a diffusion limit of the system  \eqref{eq: active leaders}-\eqref{eq:governing_eq2new}-\eqref{eq: alignment first first} and in Subsection \ref{slowsection} a diffusion limit for the system \eqref{eq: active leaders}-\eqref{eq:governing_eq2new}-\eqref{eq: alignment first first} assuming the population of leaders is ``slow" compared to the followers by scaling $c_a$, $c_p$. Finally, in Subsection \ref{sec: ranodmness passive} we develop a diffusion limit of the system \eqref{eq: active leaders}-\eqref{eq:governing_eq2new}-\eqref{eq: alignment first first} for slow leaders and including noise in the movement of passive leaders.

\subsection{Properties and asymptotic expansion of operators.} 
{ In the following sections we consider 
\begin{align*}
    T^\varepsilon=T^0+\varepsilon T^1+\mathcal{O}(\varepsilon^2)\ ,\quad \Phi^\varepsilon=\Phi^0+\varepsilon\Phi^1+\mathcal{O}(\varepsilon^2)\quad\textnormal{and}\quad B^\varepsilon=B^0+\varepsilon B^1+\mathcal{O}(\varepsilon^2)\ .
\end{align*}
Moreover, we have that
\begin{align*}
    \int_S T^0 1\ d \theta+\varepsilon\int_S T^1 1\ d\theta=1\ ,
\end{align*}
so that
\begin{align*}
    \int_S T^0 1\ d \theta = 1, \ \int_S T^1 1\ d\theta=0\ ,
\end{align*}
and similarly 
\begin{align}
    \int_S \Phi^0 d \theta=1,\ \int_S \Phi^1d\theta=0\ , \ \int_S B^0 1\ d \theta=1,\ \int_S B^11\ d\theta=0\ . \label{eq: properties of operators}
\end{align}
In Section \ref{sub: hyperbolic limit} we further show that
\begin{equation*}
    \int_S\theta\Phi^0(\Lambda\cdot\theta)d\theta=z\Lambda^0\ , \ \int_S\theta B^0d\theta=\int_S\theta B^1 d\theta=0\ .
\end{equation*}
}
\subsection{Hyperbolic limit}\label{sub: hyperbolic limit}
In this section we investigate the macroscopic dynamics of the swarm over short time scales. Consider the following scaling $(\mathbf{x},t)\mapsto(\mathbf{x}/\varepsilon,t/\varepsilon)$ where $\varepsilon=\bar{\tau}/\mathcal{T}\ll 1$ and $\mathcal{T}$ is the macroscopic time. The transition rates $R^\varepsilon_{ap}, R^\varepsilon_{pa}$ are multiplied by a factor of  $\varepsilon$ meaning that the switching rates are large enough to be observed at the macroscopic spatial scale. Then we write
\begin{align}
    \varepsilon\partial_t\rho_a-\varepsilon c_a\mathbf{b}\cdot\nabla\rho_a & =-\varepsilon R_{ap}^\varepsilon\rho_a+\varepsilon R^\varepsilon_{pa}\rho_p\label{eq: streakers hyperbolic}\ ,\\
    \varepsilon\partial_t\rho_p+\varepsilon c_p\mathbf{c}\cdot\nabla\rho_p & =\varepsilon {R}^\varepsilon_{ap}\rho_a-\varepsilon R_{pa}^\varepsilon\rho_p\label{eq: passive hyperbolic}\ ,\\
    \varepsilon\partial_t\sigma_f+\varepsilon c_f\theta\cdot\nabla\sigma_f & =-\beta\sigma_f+\zeta\beta T^\varepsilon\sigma_f+(1-\zeta)\beta\Phi^\varepsilon(\Lambda\cdot\theta)\rho_f\label{eq: followers hyperb}\ .
\end{align}

For the case of the follower population and using the expansions \eqref{eq: expansions followers} for $\sigma_f$, $T^\varepsilon$, $\Phi^\varepsilon$ and $\rho_f^\varepsilon$ from Appendix \ref{sec: Macroscopic diffusion limit} we  can rewrite (\ref{eq: followers hyperb}) as
 \begin{align}
 \varepsilon\partial_t(\sigma_f^0+\varepsilon\sigma_f^1)+\varepsilon c_f\theta\cdot\nabla(\sigma_f^0+\varepsilon\sigma_f^1)=&-\beta(\sigma_f^0+\varepsilon\sigma_f^1)+\zeta\beta(T^0+\varepsilon T^1)(\sigma_f^0+\varepsilon\sigma_f^1)\nonumber\\ &+(1-\zeta)\beta(\Phi^0+\varepsilon \Phi^1)(\rho_f^0+\varepsilon\rho_f^1)  \ .\label{eq: right hand side}
 \end{align}
For $\varepsilon^0$ we obtain
\begin{equation}
    \sigma_f^0=(\zeta+(1-\zeta)\Phi^0(\Lambda\cdot\theta))\rho_f^0\label{eq: expansion followers}\ .
\end{equation}
Here we considered that $T^0\sigma_f^0=\rho_f^0$  (see Appendix \ref{sec: turn_angle_properties} for properties of the operator $T$).

Substituting (\ref{eq: expansion followers}) into (\ref{eq: followers hyperb}) and integrating with respect to $S$ gives
\begin{align}
    \partial_t\rho_f^0\int_S(\zeta+(1-\zeta)\Phi^0(\Lambda\cdot\theta))d\theta+c_f\int_S\theta\cdot\nabla(\zeta+(1-\zeta)\Phi^0(\Lambda\cdot\theta))\rho_f^0d\theta=0 \ .
\end{align}
In the left hand side we use again $\int_S\Phi^0(\Lambda\cdot\theta)d\theta=1$ and, $\int_S\theta\Phi^0(\Lambda\cdot\theta)d\theta=z\Lambda^0$ where $z$ is given by (\ref{eq: constant z}). Hence the conservation equation for the followers reads
\begin{equation}
    \partial_t\rho_f^0+c_fz(1-\zeta)\nabla\cdot(\rho_f^0\Lambda^0)=0\ .
\end{equation}
Next we need to find the mean direction $\rho_f^0\Lambda^0$. Following the same analysis as in \cite{dimarco2016self,estrada2019interacting} we substitute the expansion for $\sigma_f^0$ given in (\ref{eq: expansion followers}) into (\ref{eq: followers hyperb}) and multiply by $\theta\cdot v$, where $v\in\mathbb{R}^n$ is orthogonal to $\Lambda^0$. Integrating over $S$ then gives
\begin{equation}
    \Bigl(\partial_t\int_S\theta\Psi(\theta)\rho_f^0d\theta+c_f\int_S\theta\cdot\nabla(\Psi(\theta)\rho_f^0)\theta d\theta\Bigr)\cdot v=\mathcal{O}(\varepsilon)\ .
\end{equation}
Here $\Psi(\theta)=\zeta+(1-\zeta)\Phi^0(\Lambda\cdot\theta)$, {$\mathcal{O}(\varepsilon)$ includes $\zeta\beta\int_S\theta T^1\sigma_f^0 d\theta$, $(1-\zeta)\beta\int_S\theta\Phi^1\rho_f^0 d\theta$ and $-\beta\int_S\sigma_f^1d\theta$}, and letting $\varepsilon\rightarrow 0$ in the right hand side we obtain
\[
\left(z(1-\zeta)\partial_t(\rho_f^0\Lambda^0)+c_f\int_S\theta\cdot\nabla(\rho_f^0\Psi(\theta))\theta d\theta \right)\cdot v=0\ .
\]

Using the fact that $v \perp \Lambda^0$, we can reformulate the above expression in terms of the orthogonal projection $P_\perp=\mathbbm{1}-\Lambda^0\otimes\Lambda^0$ onto $\Lambda^0_\perp$, 
\begin{equation}
P_\perp\Bigl( z(1-\zeta)\partial_t(\rho_f^0\Lambda^0)+c_f\nabla\cdot\rho_f^0\int_S(\theta\otimes\theta)\Psi(\theta)d\theta\Bigr)=0\ .
\end{equation}
For the first term of the above expression we can write
\begin{equation}
    z(1-\zeta)P_\perp(\rho_f^0\partial_t\Lambda^0+\Lambda^0\partial_t\rho_f^0)=z(1-\zeta)\rho_f^0\partial_t\Lambda^0\ ,
\end{equation}
since $\langle\partial_t\Lambda^0,\Lambda^0 \rangle=\frac{1}{2}\partial_t|\Lambda^0|^2=0$, i.e., $\Lambda^0\perp\partial_t\Lambda^0$. For the second term we must compute the integral $\int_S(\theta\otimes\theta)\Psi(\theta)d\theta$, where we use $\theta=\cos(s)\Lambda^0+\sin(s)\Lambda^0_\perp$ in polar coordinates for $n=2$ and spherical coordinates for $n=3$ as in \cite{dimarco2016self}. Finally, we obtain
\begin{equation}
    {\rho_f^0}(z(1-\zeta)\partial_t\Lambda^0+C_1\Lambda^0\cdot\nabla\Lambda^0)+C_2P_\perp\nabla\rho_f^0=0\ ,
\end{equation}
where  we have used $\Lambda^0_\perp\otimes\Lambda^0_\perp=\mathbbm{1}-\Lambda^0\otimes\Lambda^0$. Here $C_1=c_f(1-\zeta)a_3$ and $C_2=c_f(1-\zeta)\mathbbm{1}a_1+c_f\mathbbm{1}\pi\zeta$ for, $a_3=a_0-a_1$ and 
\begin{equation}
  a_0=\begin{cases}
    \int_{0}^{2\pi}\Phi^0(\cos(s))\cos(s)^2 ds\ , & \text{if $n=2$}\ ,\\
    2\pi\int_{0}^{\pi}\Phi^0(\cos(s))\cos(s)^2\sin(s) ds\ , & \text{if $n=3$}\ ,\label{eq: constant a0}
  \end{cases}
\end{equation}
\begin{equation}
  a_1=\begin{cases}
    \int_{0}^{2\pi}\Phi^0(\cos(s))\sin(s)^2 ds\ , & \text{if $n=2$}\ ,\\
    \pi\int_{0}^{\pi}\Phi^0(\cos(s))\sin(s)^3\sin(s) ds\ , & \text{if $n=3$}\ .\label{eq: constant a1}
  \end{cases}
\end{equation}

The final system of equations is given as
 \begin{align}
 \partial_t\rho_a^0- c_a\mathbf{b}\cdot\nabla\rho_a^0 + R_{ap}^0\rho_a^0-R^0_{pa}\rho_p^0 & =0\ ,  \label{eq: hyperbolic streakers}\\\
      \partial_t\rho_p^0+ c_p\mathbf{c}\cdot\nabla\rho_p^0 - R_{ap}^0\rho_a^0+R^0_{pa}\rho_p^0 & =0\ , \label{eq: hyperbolic passive}\\
      \partial_t\rho_f^0+c_fz(1-\zeta)\nabla\cdot(\rho_f^0\Lambda^0) &=0\ ,\label{eq: hyperbolic followers 1}\\ \rho_f^0(z(1-\zeta)\partial_t\Lambda^0 +C_1{\Lambda^0}\cdot\nabla\Lambda^0)+C_2P_\perp\nabla\rho_f^0 &=0\ .\label{eq: hyperbolic followers 2}
 \end{align}
Unlike in the parabolic limit \eqref{parstreakers}-\eqref{parfollow} in Subsection \ref{sub: diffusion limit}, we now obtain transport equations for the leader populations, which no longer adjust instantaneously to the swarm of followers.

\subsection{Hyperbolic limit for passive leaders with noise}\label{sec: passive leaders with noise}
Defining the macroscopic density 
$
\rho_p(\mathbf{x},t)=\int_0^t\int_S\bar{\sigma}_p(\cdot,\eta,\tau)d\eta d\tau
$, we study now the hyperbolic limit of \eqref{eq: passive 2} where, after integrating with respect to $\tau$ and using the boundary condition $\bar{\sigma}_p(\cdot,\eta,\tau=0)$ we write
\begin{equation}
\varepsilon\partial_t\sigma_p+\varepsilon c_p\theta\cdot\nabla\sigma_p  =\varepsilon \bar{R}^\varepsilon_{ap}(\theta)\rho_a-\varepsilon R_{pa}^\varepsilon\sigma_p-\beta\sigma_p+\beta B^\varepsilon\rho_p\ .\label{eq: hyperbolic noise}
\end{equation}

Following the same procedure as before we consider in this case  the expansions \eqref{eq: passive expansion} and \eqref{eq: other expansions} for $\sigma_p^\varepsilon$, $\rho_p^\varepsilon$, $\bar{R}_{ap}^\varepsilon$ and $R_{ap}^\varepsilon$,  and from \eqref{eq: hyperbolic noise} we obtain 
\begin{alignat}{3}
        &\varepsilon^{0}:& \ \ \sigma_p^0 &=B^0 \rho_p^0\ \label{eq: hyperbolic passive1},\\
        &\varepsilon^{1}:&\ \   -\beta\sigma_p^1  &=\partial\sigma_p^0+c_p\theta\cdot\nabla\sigma_p^0-\bar{R}_{ap}^0\rho_a^0+R^0_{pa}\sigma_p^0-\beta B^0\rho_p^1-\beta B^1\rho_p^0\nonumber\\
        &       &       & =B^0\partial_t\rho_p^0+c_pB^0\theta\cdot\nabla\rho_p^0-\bar{R}_{ap}^0\rho_a^0+R^0_{pa}B^0\rho_p^0\nonumber\\ &     &   &\ \ \ \ \ \ \ \ \ -\beta B^0\rho_p^1-\beta B^1{\rho_p^0}\ . \label{eq: hyperbolic passive2}
\end{alignat}
The next step is to compute the term $\partial_t\rho_p^0$. For that, we integrate \eqref{eq: hyperbolic noise} with respect to $\theta$ and we have
\begin{equation}
\varepsilon\partial_t(\rho_p^0+\varepsilon\rho_p^1)+{\varepsilon c_p}\nabla\cdot\int_S\theta(\sigma_p^0+\varepsilon\sigma_p^1)d\theta= \varepsilon(R^0_{ap}+\varepsilon R_{ap}^1)(\rho_a^0+\varepsilon\rho_a^1)-\varepsilon(R^0_{pa}+\varepsilon R^1_{pa})(\rho_p^0+\varepsilon\rho_p^1)\label{eq: conservation leaders hyperbolic}\ .
\end{equation}
{ Here we have used \eqref{eq: properties of operators}. With \eqref{eq: conservation leaders hyperbolic}} and since $\int_S\theta B^0d\theta=0$,  we obtain
\[
\partial_t\rho_p^0=R_{ap}^0\rho_a^0-R^0_{pa}\rho_p^0\ .
\]
Substituting the above expression in \eqref{eq: hyperbolic passive2} we finally obtain
\begin{equation}
\sigma_p^1=\frac{-1}{\beta}\Bigl(B^0R^0_{ap}\rho_a^0+c_pB^0\theta\cdot\nabla\rho_p^0-{\bar{R}_{ap}^0}\rho_a-\beta B^0\rho_p^1-\beta B^1\rho_p^0 \Bigr)\ ,\label{eq: expansion for sigma_1}
\end{equation}
and hence, $\sigma_p=B^0\rho_p^0+\varepsilon\sigma_p^1+\mathcal{O}(\varepsilon^2)$. Going back to the conservation equation \eqref{eq: conservation leaders hyperbolic} we can write now
\begin{equation}
    \partial_t\rho_p^0-R^0_{ap}\rho_a^0+R^0_{pa}\rho_p^0=0\ .
\end{equation}
 Computing the term 
 $
{-\varepsilon c_p}\nabla\cdot\int_S\theta\sigma_p^1d\theta
 $
 using \eqref{eq: expansion for sigma_1} we obtain that
\begin{equation}
    {-\varepsilon c_p}\nabla\cdot\int_S\theta\sigma_p^1d\theta=\varepsilon (\Delta(D_p\rho_p^0)-\nabla\cdot(\mathcal{R}_{ap}^0\rho_a^0))\ ,
\end{equation}
where $D_p=\frac{c_p}{\beta}\int_S\theta\theta^TB^0d\theta$ and $\mathcal{R}_{ap}^0=\frac{c_p}{\beta}\int_S\theta\bar{R}_{ap}^0d\theta$, which shows that the diffusion term only appears for higher orders of $\varepsilon$. The hyperbolic system of equations reads as
 \begin{align}\label{hyperbolicstreakerswithnoise}
 \partial_t\rho_a^0- c_a\mathbf{b}\cdot\nabla\rho_a^0 + R_{ap}^0\rho_a^0-R^0_{pa}\rho_p^0 & =0\ ,  \\\
     \partial_t\rho_p^0-R_{ap}^0\rho_a^0+R^0_{pa}\rho_p^0&=0\ , \label{hyperbolicpassivewithnoise}\\
      \partial_t\rho_f^0+c_fz(1-\zeta)\nabla\cdot(\rho_f^0\Lambda^0) &=0\ ,\\ \rho_f^0(z(1-\zeta)\partial_t\Lambda^0 +C_1\Lambda\cdot\nabla\Lambda^0)+C_2P_\perp\nabla\rho_f^0 &=0\ .\label{hyperbolicfollowerswithnoise}
 \end{align}
\\

\subsection{Diffusion limit}\label{sub: diffusion limit}

Next we turn to determining macroscopic equations that describe the interactions between the three different populations over long time regimes. We introduce a parabolic scaling $(\mathbf{x},t)\mapsto (\mathbf{x}/\varepsilon,\ t/\varepsilon^2 )$, where $\varepsilon$ is defined in the sense of Section \ref{sub: hyperbolic limit}.  

Equations (\ref{eq: streakers}), (\ref{eq: passive}) and (\ref{eq: followers}) then become
\begin{align}
    \varepsilon^2\partial_t\rho_a-\varepsilon c_a\mathbf{b}\cdot\nabla\rho_a & =-\varepsilon R_{ap}^\varepsilon\rho_a+\varepsilon R_{pa}^\varepsilon\rho_p\ \label{eq: streakers original},\\
  \varepsilon^2\partial_t\rho_p+\varepsilon c_p\mathbf{c}\cdot\nabla\rho_p & =\varepsilon {R}_{ap}^\varepsilon\rho_a-\varepsilon R_{pa}^\varepsilon\rho_p\label{eq: passive scaled} \ ,\\
    \varepsilon^2\partial_t\sigma_f+\varepsilon c_f\theta\cdot\nabla\sigma_f & =-\beta\sigma_f+\zeta\beta T^\varepsilon\sigma_f+\varepsilon(1-\zeta)\beta\Phi^\varepsilon(\Lambda\cdot\theta)\rho_f\ .\label{eq: diffusion followers}
\end{align}
{Note that we have also scaled the alignment kernel $\Phi^\varepsilon(\Lambda\cdot\theta)$.}

The conservation equation (\ref{eq: conservation followers}) for the follower population is
\begin{equation}
\partial_t\rho_f+nc_f\nabla\cdot w_f=0\ ,\label{eq: conservation equation followers}
\end{equation}
where $w_f=\frac{1}{n}\int_S\theta\sigma_fd\theta$. The next step is to compute the mean direction $w_f$. We follow the steps in Appendix \ref{sec: Macroscopic diffusion limit} and, using asymptotic expansions for  the terms $\sigma_f$, $\rho_f$, $T^\varepsilon$ and $\Phi^\varepsilon$ we can re-write the conservation equation \eqref{eq: conservation equation followers} as 
\begin{equation*}
    \varepsilon^2\partial_t\rho_f^0+{\varepsilon^2c_f}\nabla\cdot\int_S\theta\sigma_f^1d \theta=0\ .
\end{equation*}
Finally, we obtain
\begin{equation}
\partial_t\rho_f^0-c_f\Delta(D_f\rho_f^0)+\nabla\cdot(z\Lambda^0_W\rho^0_f)=0\ ,\label{eq: followers diffusion final}
\end{equation}
where $z\Lambda^0_W=\int_S\theta\Phi^0({\Lambda_W^0}\cdot\theta)d \theta$, $D_f=\frac{\zeta c_f}{\beta(1-\zeta)}\int_S\theta\theta^Td\theta$, $z$ is given by \eqref{eq: constant z} and $W$ is the total mean direction of the whole population.

The system describing the macroscopic densities of followers and leaders, in the diffusion limit, reads as follows
\begin{align}
-c_a\mathbf{b}\cdot\nabla\rho_a^0+ R_{ap}^0\rho_a^0- R_{pa}^0\rho_p^0 & =0\ , \label{parstreakers}\\
     c_p\mathbf{c}\cdot\nabla\rho_p^0- R_{ap}^0\rho_a^0+ R_{pa}^0\rho_p^0 & =0\ ,  \label{parpassive}\\ \partial_t\rho_f^0+c_f\nabla\cdot w_f&=0\\
     w_f=-\nabla D_f\rho_f^0&+\rho_f^0z\Lambda^0_W\ .  \label{parfollow}
\end{align}
Note that the evolution of the swarm is determined by the movement of the followers. There are no time derivatives in the equations for the densities of the leaders, which adjust instantaneously to the movement of the followers. Technically, the time derivatives in the system \eqref{eq: streakers original}, \eqref{eq: passive scaled}, \eqref{eq: diffusion followers} scale as $\varepsilon^2$, so that the leader's movement is determined at order $\varepsilon$ by convection and the transition between their active and passive states.

\subsection{Diffusion limit for slow leaders}\label{slowsection}

We consider the situation where $c_a$, $c_p$, $R_{ap}^\varepsilon$ and $R_{pa}^\varepsilon$ are all of order $\varepsilon$, in particular the leaders move slowly compared to the followers.

Equations (\ref{eq: streakers}), (\ref{eq: passive}) and (\ref{eq: followers}) then become
\begin{align}
    \varepsilon^2\partial_t\rho_a-\varepsilon^2 c_a\mathbf{b}\cdot\nabla\rho_a & =-\varepsilon^2 R_{ap}^\varepsilon\rho_a+\varepsilon^2 R_{pa}^\varepsilon\rho_p\ \label{eq: streakers original2},\\
  \varepsilon^2\partial_t\rho_p+\varepsilon^2 n c_p\mathbf{c}\cdot\nabla\rho_p & =\varepsilon^2 {R}_{ap}^\varepsilon\rho_a-\varepsilon^2 R_{pa}^\varepsilon\rho_p\label{eq: passive scaled slow leaders} \ ,\\
    \varepsilon^2\partial_t\sigma_f+\varepsilon c_f\theta\cdot\nabla\sigma_f & =-\beta\sigma_f+\zeta\beta T^\varepsilon\sigma_f+(1-\zeta)\beta\Phi^\varepsilon(\Lambda\cdot\theta)\rho_f\ .\label{eq: diffusion followers2}
\end{align}

The system describing the macroscopic densities of followers and leaders, in the diffusion limit, reads as follows
\begin{align}
\partial_t\rho_a^0-c_a\mathbf{b}\cdot\nabla\rho_a^0+ R_{ap}^0\rho_a^0- R_{pa}^0\rho_p^0 & =0\ , \label{parstreakersslow}\\
\partial_t\rho_p^0 +    c_p\mathbf{c}\cdot\nabla\rho_p^0- R_{ap}^0\rho_a^0+ R_{pa}^0\rho_p^0 & =0\ ,  \label{parpassiveslow}\\ \partial_t\rho_f^0+c_f\nabla\cdot w_f&=0\\
     w_f=-\nabla D_f\rho_f^0&+\rho_f^0z\Lambda^0_W\ .  \label{parfollowslow}
\end{align}

\subsection{Random movement for passive leaders}\label{sec: ranodmness passive}

We start from the mesoscopic equation \eqref{eq: hyperbolic noise} obtained from \eqref{eq: passive 2}.
Introducing the diffusive scaling as in Subsection \ref{sub: diffusion limit} we get
\begin{align}
    \varepsilon^2\partial_t\rho_a-\varepsilon^2 c_a\mathbf{b}\cdot\nabla\rho_a & =-\varepsilon^2 R_{ap}^\varepsilon\rho_a+\varepsilon^2 R_{pa}^\varepsilon\rho_p\ ,\\
  \varepsilon^2\partial_t\sigma_p+\varepsilon c_p\theta\cdot\nabla\sigma_p  &=\varepsilon^2 \bar{R}_{ap}^\varepsilon(\theta)\rho_a-\varepsilon^2 R_{pa}^\varepsilon\sigma_p-\beta\sigma_p+\beta B^\varepsilon\rho_p\label{eq: passive scaled2}\ ,\\
    \varepsilon^2\partial_t\sigma_f+\varepsilon c_f\theta\cdot\nabla\sigma_f & =-\beta\sigma_f+\zeta\beta T^\varepsilon\sigma_f+\varepsilon(1-\zeta)\beta\Phi^\varepsilon(\Lambda\cdot\theta)\rho_f\ .
\end{align}
Note that the speed $c_a$, as well as the switching rates $R^\varepsilon_{ap},\ R^\varepsilon_{pa}$ and $\bar{R}_{ap}^\varepsilon(\theta)$ are scaled as $\varepsilon^2$. The conservation equation for the passive leaders \eqref{eq: passive scaled2}, is then
\begin{equation}
    \varepsilon^2\partial_t\rho_p+{\varepsilon^2  c_p}\nabla\cdot \int_S\theta\sigma_pd\theta=\varepsilon^2R^\varepsilon_{ap}\rho_a-\varepsilon^2R_{pa}^\varepsilon\rho_p \label{eq: final conservation}\ ,
\end{equation}
where
\begin{equation}
    \rho_p=\int_S\sigma_pd\theta\ \ \textnormal{and}\ \ w_p=\frac{1}{n}\int_S\theta\sigma_pd\theta\ .\label{eq: macroscopic}
\end{equation}
Following the steps described in Appendix \ref{sec: Diffusion limit for passive leaders} we find an asymptotic expansion for $\sigma_p$, using \eqref{eq: epsilon-1} and \eqref{eq: epsilon0}, as follows
$$
\sigma_p=B^0\rho_p^0+\varepsilon\sigma_p^1+\mathcal{O}(\varepsilon^2)\ ,
$$
where $\sigma_p^1$ is given by \eqref{eq: epsilon0}.  Computing the mean direction term $\nabla\cdot\int_S\theta\sigma_p^1d\theta$ as in \eqref{eq: mean direction 2} we finally obtain, for the passive leaders with random movement,
\begin{equation}
    \partial_t\rho_p^0-c_p\Delta(D_p{+\mathbb{T}})\rho^0_p=R^0_{ap}\rho_a^0-R^0_{pa}\rho_p^0\ ,
\end{equation}
where $$D_p=\frac{pc_p}{\beta}\int_S\theta\theta^TT^0(\theta) d\theta\ ,\ \ \textnormal{and}\ \  {\mathbb{T}=\frac{(1-p)c_p}{\beta}\int_S\theta \theta^T T^0_{\rho_f}d\theta}\ .$$
Hence the final system is given by
\begin{align}
{\partial_t\rho_a^0}-c_a\mathbf{b}\cdot\nabla\rho_a^0+ R_{ap}^0\rho_a^0- R_{pa}^0\rho_p^0 & =0\ ,\label{parstreak2}\\ \partial_t\rho_p^0-c_p\Delta(D_p{+\mathbb{T}})\rho^0_p-R^0_{ap}\rho_a^0+R^0_{pa}\rho_p^0&=0\ , \label{parpassive2}\\
\partial_t\rho_f^0+c_f\nabla\cdot w_f&=0\\
     w_f=-\nabla D_f\rho_f^0&+\rho_f^0z\Lambda^0_W\ .  \label{parfollowers2}
\end{align}

\section{Numerical solution}

As an example of the macroscopic dynamics we illustrate the time evolution of the swarm dynamics for the system \eqref{parstreak2} - \eqref{parfollowers2} in a one dimensional domain $x \in (0,L)$, with Neumann conditions imposed at the boundary. The system is solved numerically using a simple finite difference scheme on a sufficiently fine equidistant mesh. More precisely, for the leaders we couple a Lax-Wendroff discretisation for the transport equation \eqref{parstreak2} to a Crank-Nicolson time stepping for the diffusion equations \eqref{parpassive2}, \eqref{parfollowers2} in a way that conserves the total number of leaders.   

Initially we co-localise leaders and followers in $(0,L) = (0,15)$, setting initial distributions proportional to a Gaussian of width $0.6$ and centre $x=3$, with a total density of $3$ for active and passive leaders, respectively $100$ for the followers. This reflects certain real world follower-leader systems, where leaders form a fairly small fraction (e.g. \cite{seeley}). The time evolution of the swarm is illustrated for model parameters $c_a \textbf{b}= 6$, $c_p (D_p{+\mathbb{T}}) = c_f D_f =1$ and $\Lambda_W^0$ as in \eqref{eq: mean direction alt form} which is determined by the leader population alone (formally, $\lambda \to \infty$). The conversion rates between active and passive leaders are taken to be 
\[
R_{ap}^0 = R^0\left(\exp(-\tfrac{\rho^0_f}{\Theta\max_{x \in (0,L)}\rho^0_f}) - \exp(-1/\Theta)\right)
\]
when $\partial_x \rho^0_f >0$ and $R_{ap}^0 = 0$ otherwise, respectively, 
\[
R_{pa}^0 = R^0 \left(\exp(-\tfrac{\rho^0_f}{\Theta\max_{x \in (0,L)}\rho^0_f}) - \exp(-1/\Theta)\right)
\]
when $\partial_x \rho^0_f <0$ and $R_{pa}^0 = 0$ otherwise. For illustration we use a cut-off fraction of $\Theta = 0.2$ and strength $R^0=75$.

\begin{figure}
\begin{subfigure}{.3\textwidth}
  \centering
  \includegraphics[width=\linewidth]{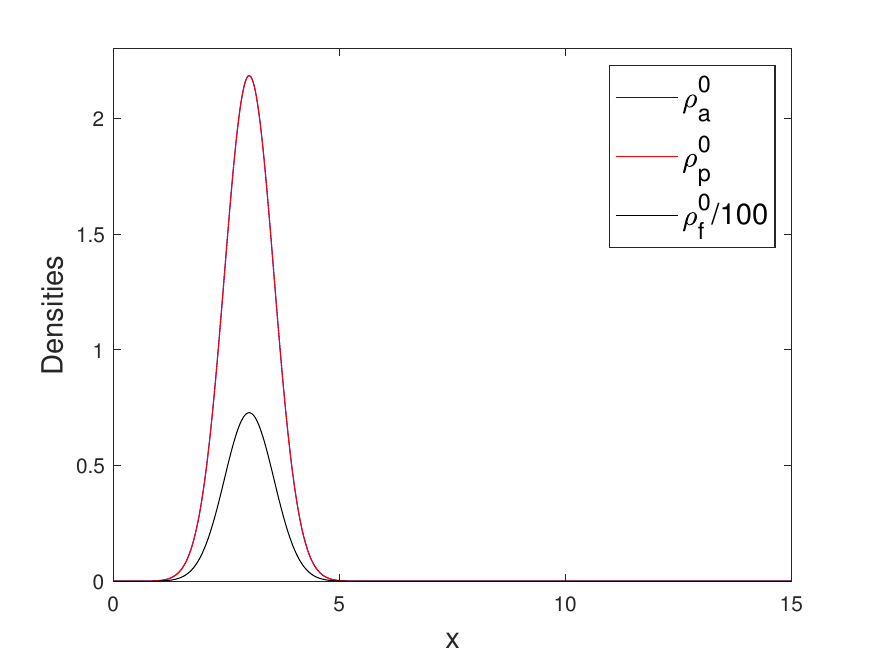}  
  \caption{$t=0$}
\end{subfigure}
\begin{subfigure}{.3\textwidth}
  \centering
  \includegraphics[width=\linewidth]{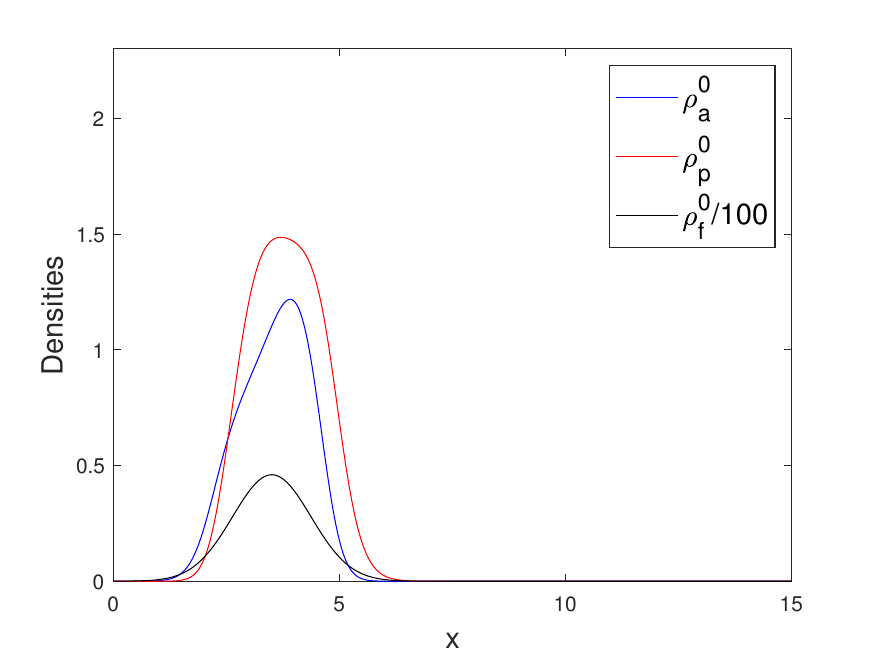}  
  \caption{$t=0.25$}
\end{subfigure}
\begin{subfigure}{.3\textwidth}
  \centering
  \includegraphics[width=\linewidth]{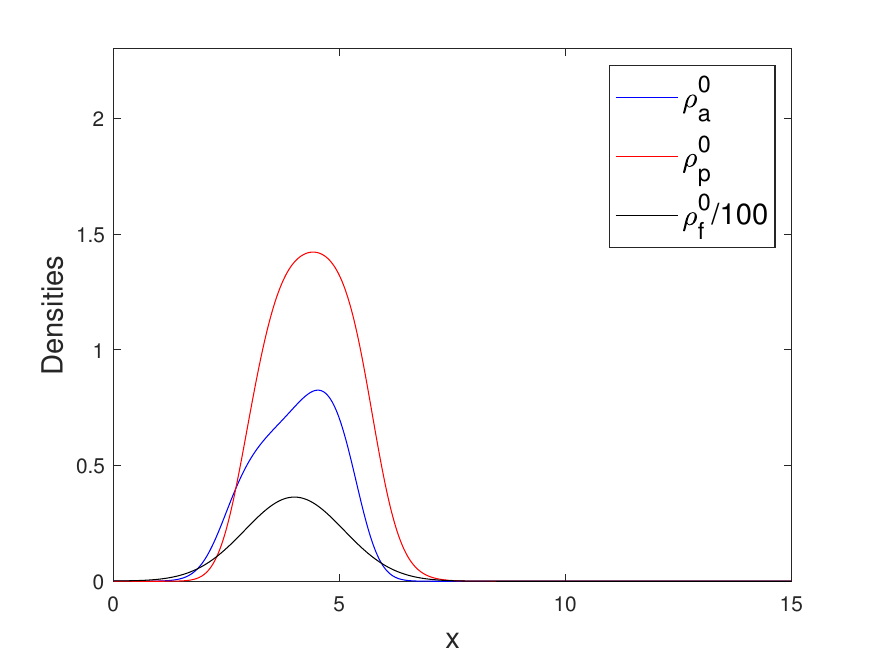}  
  \caption{$t=0.5$}
\end{subfigure}

\begin{subfigure}{.3\textwidth}
  \centering
  \includegraphics[width=\linewidth]{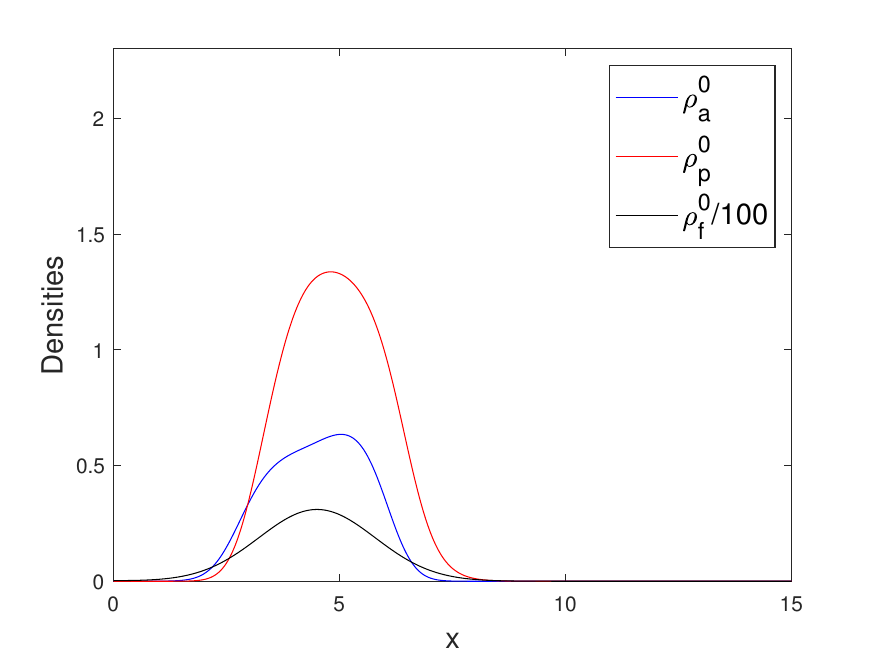}  
  \caption{$t=0.75$}
\end{subfigure}
\begin{subfigure}{.3\textwidth}
  \centering
  \includegraphics[width=\linewidth]{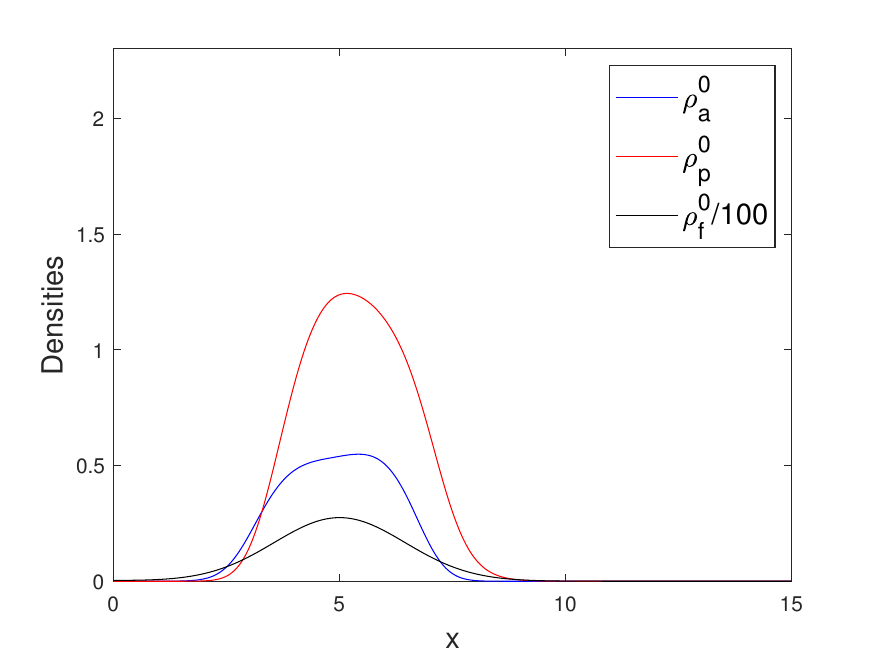}  
  \caption{$t=1$}
\end{subfigure}
\begin{subfigure}{.3\textwidth}
  \centering
  \includegraphics[width=\linewidth]{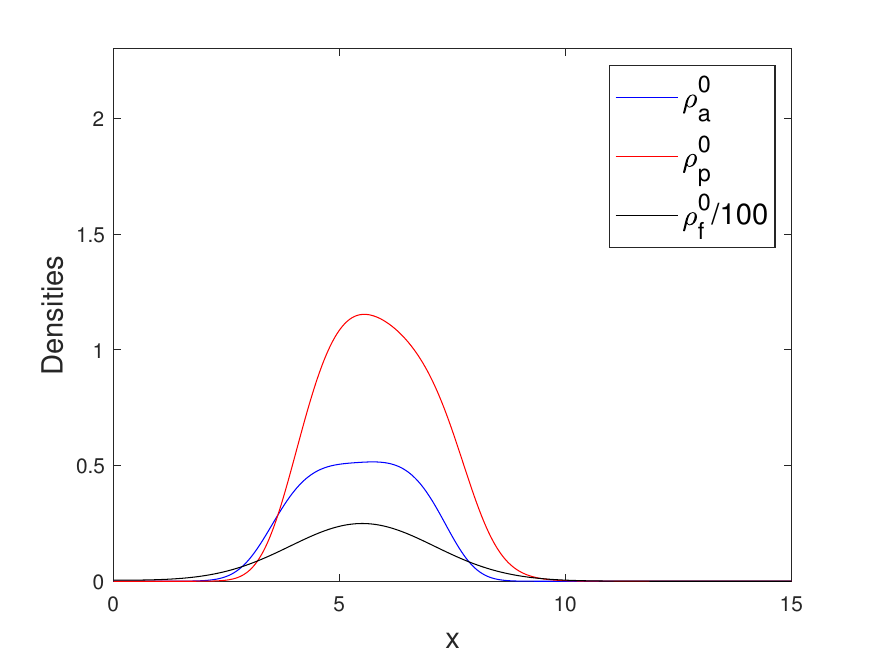}  
  \caption{$t=1.25$}
\end{subfigure}

\begin{subfigure}{.3\textwidth}
  \centering
  \includegraphics[width=\linewidth]{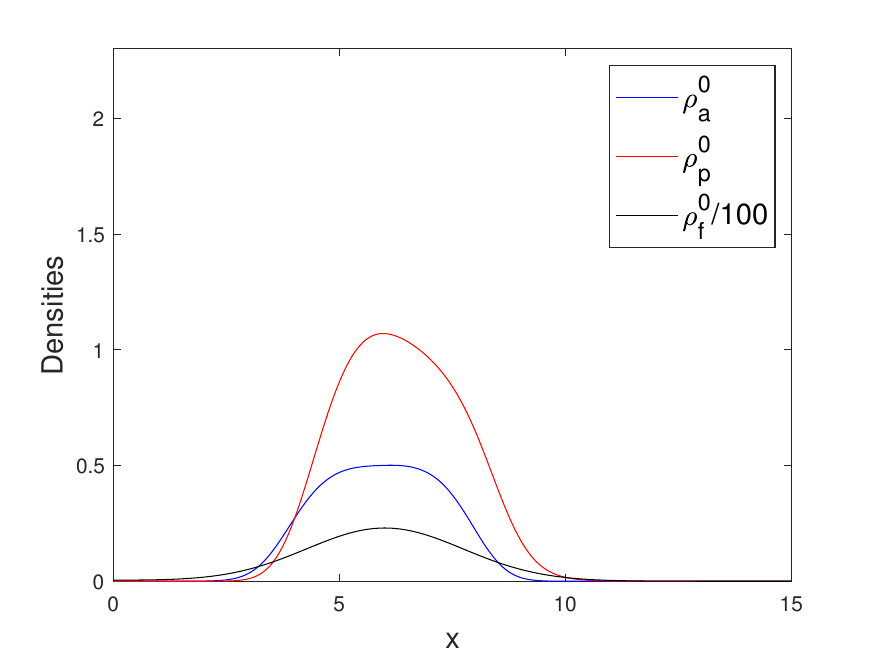}  
  \caption{$t=1.5$}
\end{subfigure}
\begin{subfigure}{.3\textwidth}
  \centering
  \includegraphics[width=\linewidth]{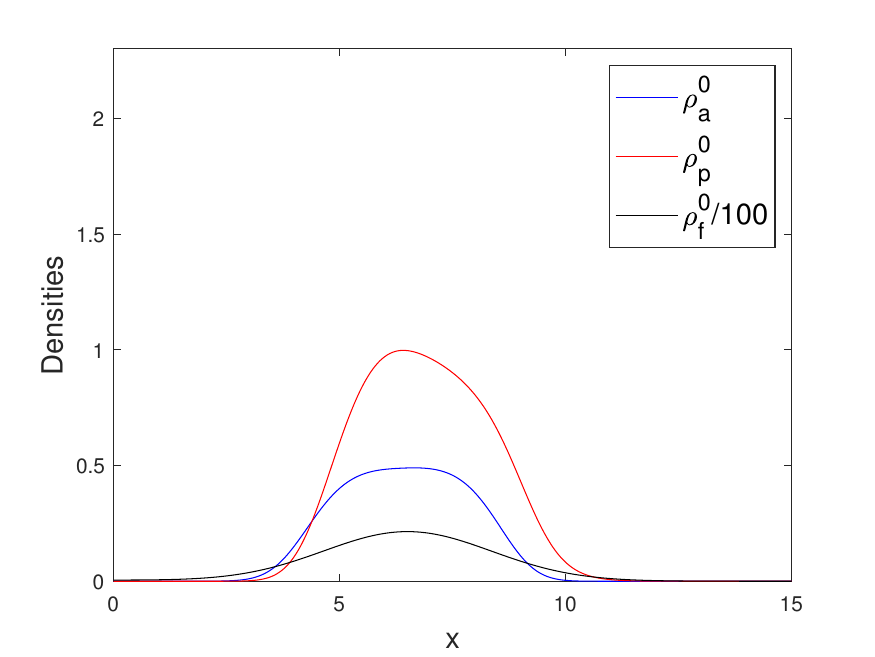}  
  \caption{$t=1.75$}
\end{subfigure}
\begin{subfigure}{.3\textwidth}
  \centering
  \includegraphics[width=\linewidth]{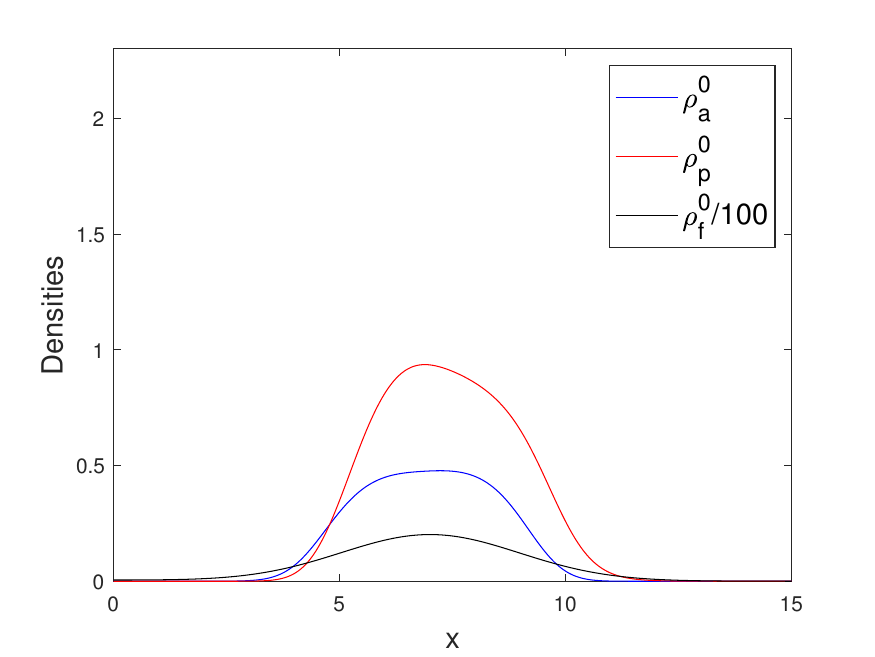}  
  \caption{$t=2$}
\end{subfigure}

\begin{subfigure}{.3\textwidth}
  \centering
  \includegraphics[width=\linewidth]{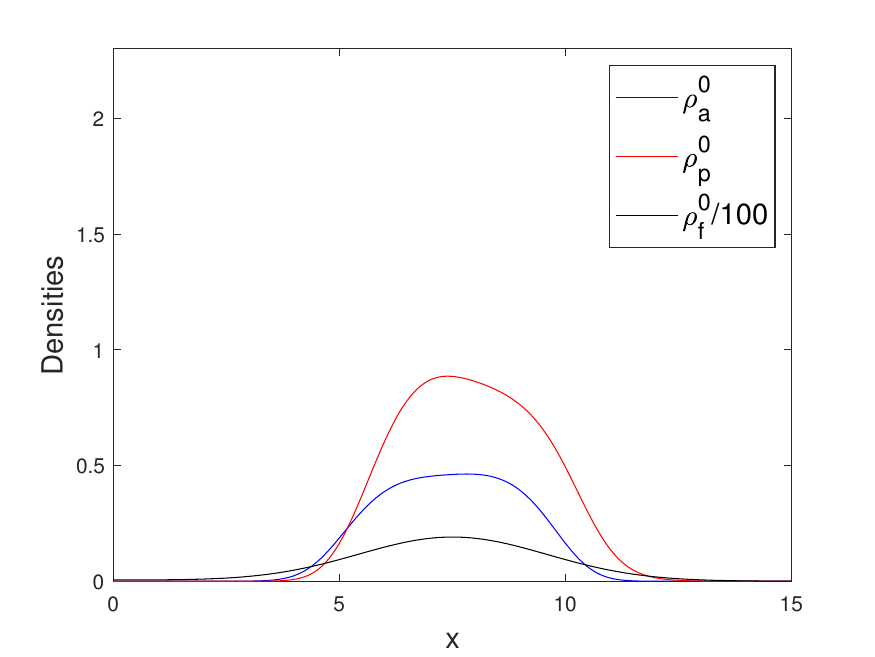}  
  \caption{$t=2.25$}
\end{subfigure}
\begin{subfigure}{.3\textwidth}
  \centering
  \includegraphics[width=\linewidth]{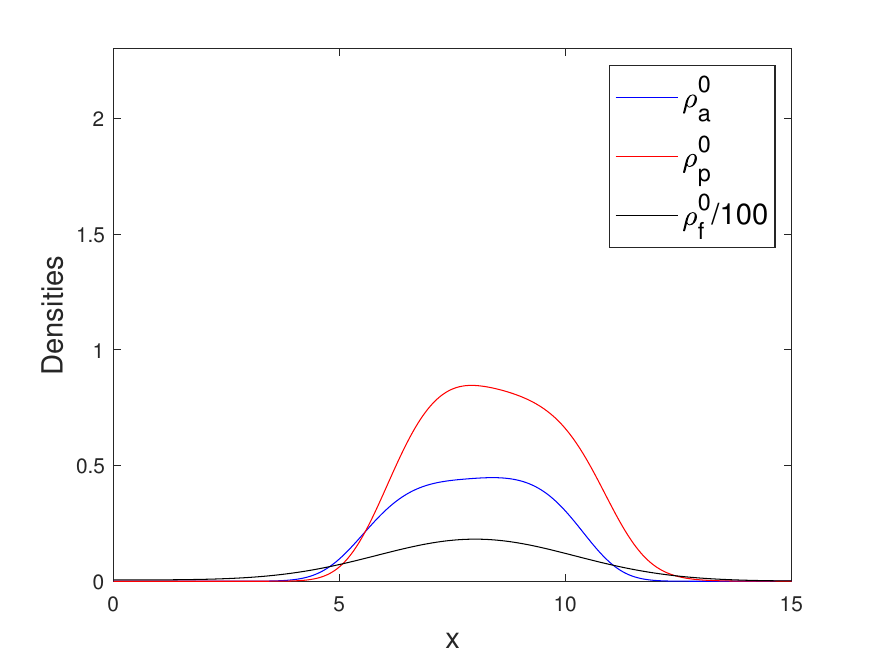}  
  \caption{$t=2.5$}
\end{subfigure}
\begin{subfigure}{.3\textwidth}
  \centering
  \includegraphics[width=\linewidth]{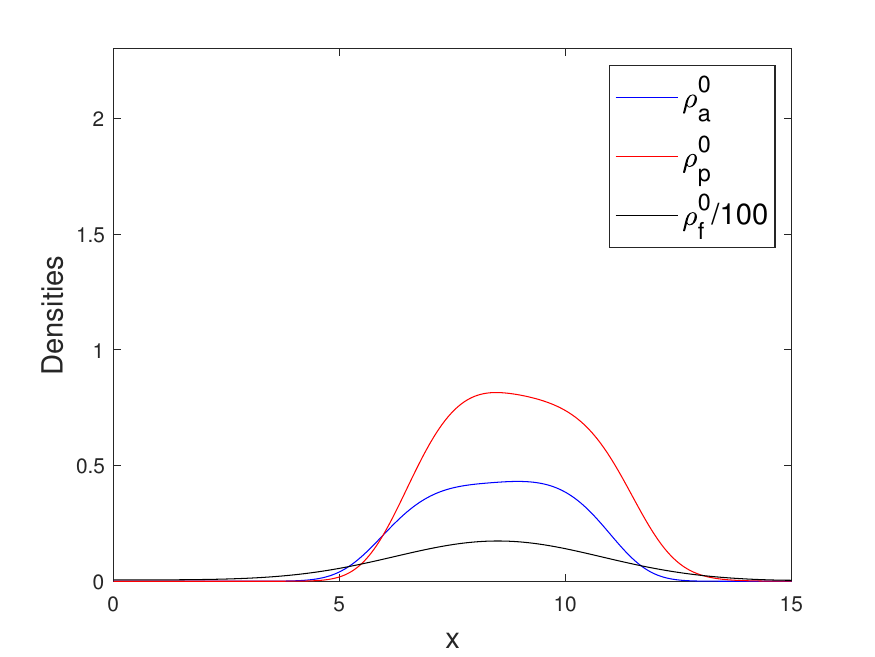}  
  \caption{$t=2.75$}
\end{subfigure}

\begin{subfigure}{.3\textwidth}
  \centering
  \includegraphics[width=\linewidth]{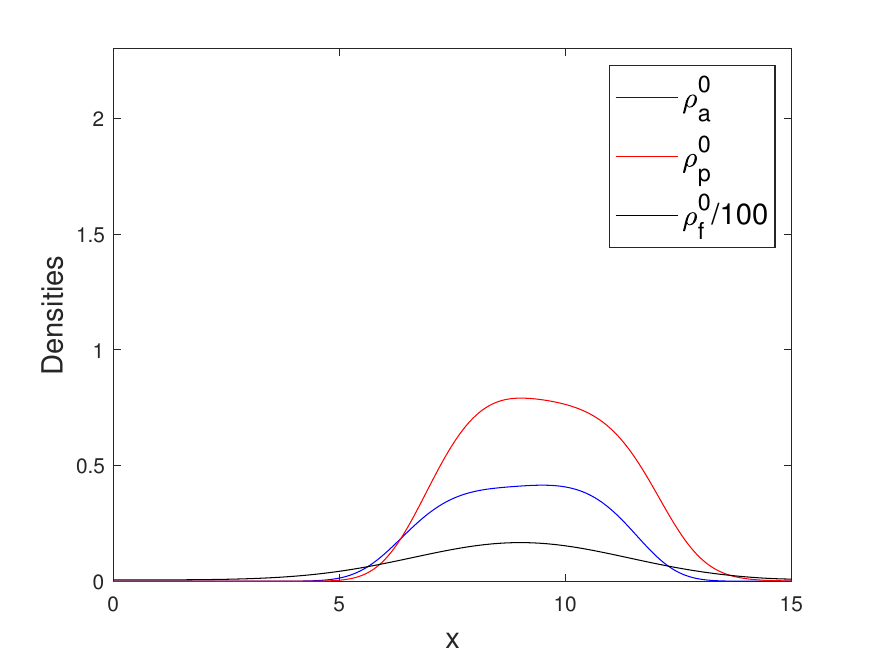}  
  \caption{$t=3$}
\end{subfigure}
\begin{subfigure}{.3\textwidth}
  \centering
  \includegraphics[width=\linewidth]{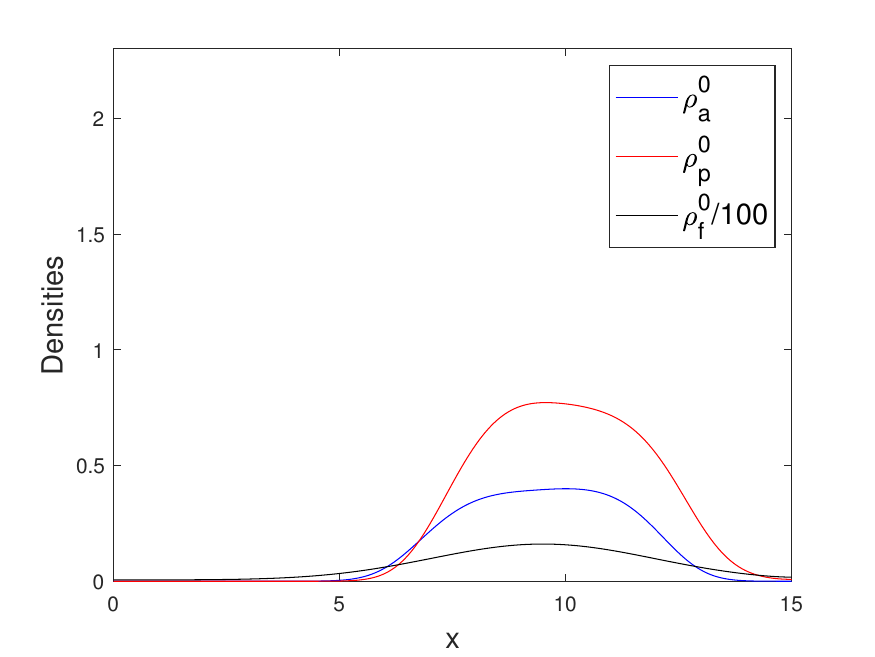}  
  \caption{$t=3.25$}
\end{subfigure}
\begin{subfigure}{.3\textwidth}
  \centering
  \includegraphics[width=\linewidth]{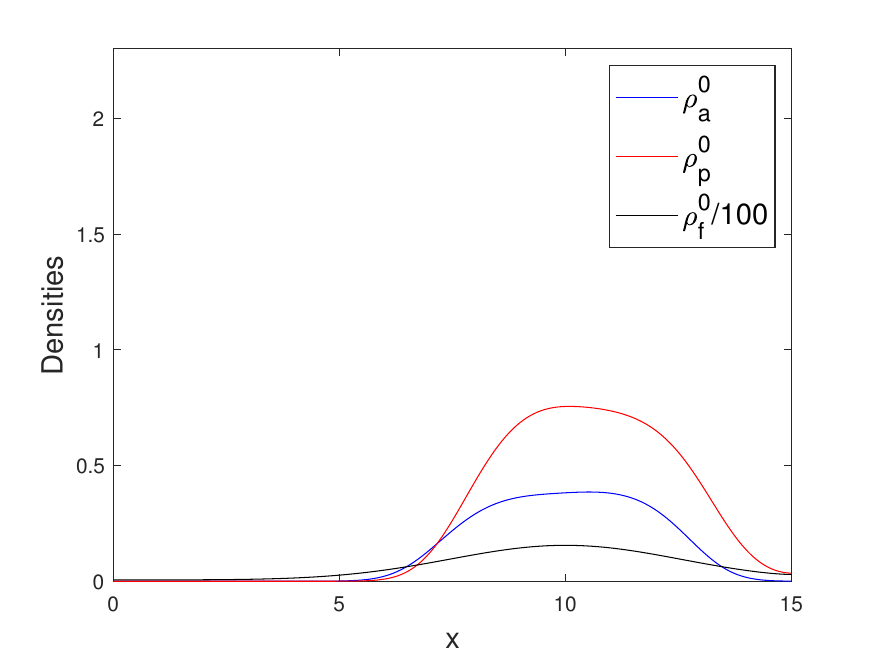}  
  \caption{$t=3.5$}
\end{subfigure}
\caption{Evolution of follower and leader populations in model example.}\label{numericalfig}
\end{figure}

\begin{figure}
  \centering
  \includegraphics[width=0.6\linewidth]{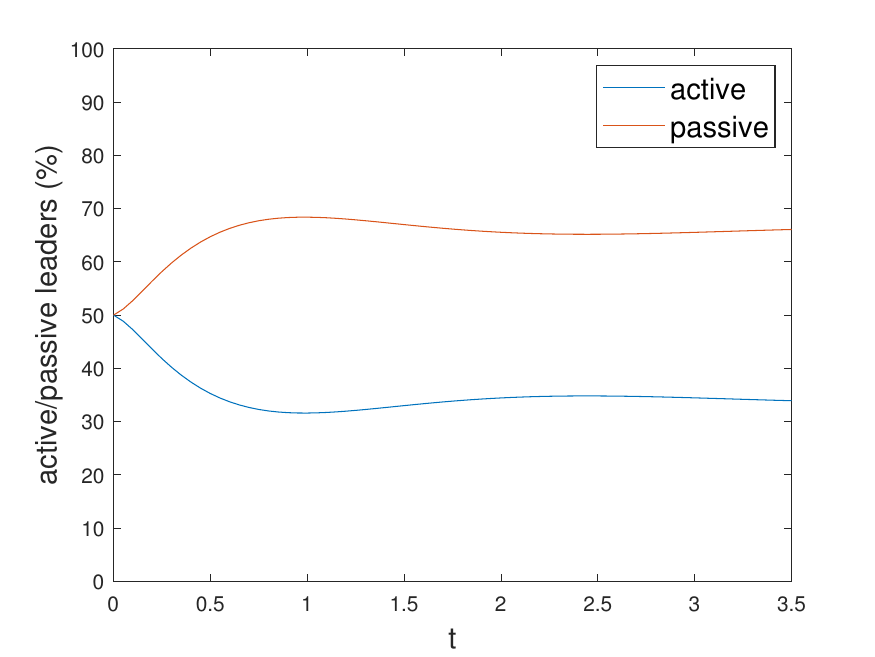}  

\caption{Percentage of active and passive leaders as a function of time.}\label{numericalfig2}
\end{figure}

Figure \ref{numericalfig} shows the resulting evolution of the densities of followers (black), active (blue) and passive (red) leaders. For purposes of illustration the density $\rho^0_f$ of followers is divided by a factor $100$.  The follower and leader populations initially move together as a localised swarm. In the current model the follower population becomes less localised with time, because the diffusion in the current model \eqref{parfollowers2} is not compensated by localising mechanisms. The leader populations widen correspondingly. Note that the fraction of active leaders quickly decreases to a small, but stable fraction, as passive leaders are slow to diffuse to the back of the swarm, where they become active again. This is confirmed by the corresponding percentages of active, respectively passive leaders shown as a function of time in Figure \ref{numericalfig2}. Because of the steep gradient in the conversion rates $R_{pa}^0$ and $R_{ap}^0$, active and passive leaders remain in the centre of the swarm, with passive leaders generated at the front. 

\section{Discussion of macroscopic equations} \label{sec:macrodiscussion}

The macroscopic equations in the hyperbolic limit, \eqref{eq: hyperbolic streakers} to  \eqref{eq: hyperbolic followers 2},  capture the behaviours depicted in Figures \ref{fig: swarm1} and 
\ref{fig: swarm23} on hyperbolic time scales.  Active leaders follow a transport equation \eqref{eq: hyperbolic streakers}, while passive leaders are transported in the opposite direction to that of the active leaders with its respective velocity. At the front of the swarm, the lower bound $R^{0}_{ap}\geq r_0$ forces $\rho_a$ to decay exponentially fast to zero over a characteristic length scale  $\frac{c_a}{r_0}$, as active leaders convert into passive leaders. As passive leaders reach areas of large $R^{0}_{pa}$, corresponding to the rear of the swarm, they convert into active and head once again to the swarm front. Here the characteristic length scale is $\frac{c_p}{r_0}$. Note that the movement of active and passive leaders depends on followers only through the transition rates $R_{ap}^0$ and $R_{pa}^0$: this is logical enough, given that leaders have knowledge of the target site and should not be swayed by the uninformed. A more detailed model may also incorporate some influence of the followers on the leader direction (e.g. due to avoiding collisions), however we have neglected that here for simplicity.

The equations describing dynamics of the followers,  \eqref{eq: hyperbolic followers 1} and \eqref{eq: hyperbolic followers 2}, correspond to classical equations for swarming particles \cite{estrada2019interacting}. The directional movement of the swarm depends on leaders only through $\Lambda^W$, and the swarm impacts on the movement of the leaders only through the switching rates $R_{ap}^0$ and $R_{pa}^0$ as noted above. This illustrates a limitation of the homogeneous choice for the alignment kernel \eqref{eq: mean direction}: since 
$\Lambda^W$ is a unit vector, only a single leader is required to direct the swarm to its target, a result which clearly stretches credulity for very large swarms. To avoid such pathological cases it is possible to consider inhomogeneous alignment kernels, see Section \ref{sec: inhom}, where the size of the orientation vector  $\Lambda^W$ is taken to increase both with the number of leaders moving along $\mathbf{b}$ and with the strength of the interactions.

A fundamental question concerns whether two follower populations, as in Figure \ref{fig: swarm23}, behave as one joint swarm or two separate swarms. Building on the discussion introduced in Section \ref{sec: swarm description}, for the model here the distinction lies with the separation of the two populations: for separations much larger than $\frac{c_a}{r_0}$, leaders remain confined to their separate follower swarms; for 
separations smaller than $\frac{c_a}{r_0}$, the leader populations bridge the gap that separates the peaks to join the swarms together; active to passive leader conversion occurs predominantly within a single zone at the swarm front.

Equations \eqref{hyperbolicstreakerswithnoise} to  \eqref{hyperbolicfollowerswithnoise}, where a heading noise is incorporated for passive leaders, support the general description above with a few noteworthy modifications. The passive leaders now undergo an additional undirected diffusion through the swarm, and hence their dynamics follow an advection-diffusion equation \eqref{hyperbolicpassivewithnoise} rather than the transport equation above.  The terms of order $\varepsilon$ account for secondary effects in the collective movement. 


Describing the swarm on parabolic time scales, however, proves more intricate. The macroscopic equations indicate some inherent problems in describing both the advective movement of active leaders and the diffusive movement of the followers. Unlike the hyperbolic limits, the density of followers in \eqref{parfollow} is now observed to diffuse over parabolic time scales, implying an eventual loss of swarm cohesion. Further, equations \eqref{parstreakers} and \eqref{parpassive} may not conserve the total number of leaders $\int (\rho_a^0 + \rho_p^0)\ dx$.

Nevertheless, the general shape of the leader population exhibits the form expected for a swarm. To see this, first note that the distribution of active leaders on long (parabolic) time scales is constant along the vector field $\mathbf{b}$ whenever $R^{0}_{ap}$ and $R^{0}_{pa}$ vanish, i.e.~within the interior of the swarm. Further, even within regions of nonzero $R^{0}_{ap}$ and $R^{0}_{pa}$ we note that 
$$-c_a\mathbf{b}\cdot\nabla\rho_a^0 + c_p\mathbf{c}\cdot\nabla\rho_p^0 = 0 \ .$$
For example, in one dimension this implies that  $\rho_p^0 = \frac{c_a}{c_p}\rho_a^0$, independent of the transition rates $R^{0}_{ap}$ and $R^{0}_{pa}$. Therefore
$$c_a\partial_x \rho_a^0 = (R^{0}_{ap} - \frac{c_a}{c_p} R^{0}_{pa})\rho_a^0\ .$$
The exponential decay of the density in front of the swarm and behind the swarm easily follows. 

When leader populations are assumed to move slowly we derive the parabolic model \eqref{parstreakersslow} to \eqref{parfollowslow}. Unlike the first parabolic scaling attempt, this system conserves the total number of leaders and the leader populations and the dynamics follow the transport equations known for hyperbolic time scales. {Swarms successfully guided through slow leaders have a logical basis, for example consider some adults attempting to guide a group of small children through a park: the latter tend to move in a rather chaotic and random way, so the adults are forced to slow their movements to keep the group intact. The assumption of slow leaders does not hold, however, for certain other swarming systems, an example being bee swarms where active leaders are believed to engage in fast streaking behaviour.} Fast directed movements prove problematic to handle in the parabolic limit, which can be attributed to the disparate time scales involved: fast target-directed streaks would result in an active leader rapidly covering a swarm's dimension, hence resulting (on a parabolic time scale) on almost instantaneous transitions from passive to active and back to passive leader. The slow leader model, can perhaps be regarded as a characterisation of such systems in which the circular movement of the leaders (moving to the front of the swarm and back) played out over a number of cycles is regarded as an effective translational movement of the leader population, now slow compared to the random movement of the followers. A careful analysis might, by separating fast and slow scales of the movement, be reduced to this situation. Such an analysis is beyond the scope of the current work.

Similar conclusions are obtained for the model in Section \ref{sec: ranodmness passive}, corresponding to slow streakers and passive leaders which follow a random walk. Here, the transport equation for the streakers is combined with a diffusion equation for the passive leaders. The  total number of leaders is conserved.

The modelling difficulties in the parabolic limit may be circumvented by \emph{prescribing} the macroscopic evolution of the leader population, so that all properties macroscopically expected for the leader population are imposed in a phenomenological model. It replaces Equations (\ref{eq: streakers}), (\ref{eq: passive})  by a deterministic set of equations for $\rho_a$ and $\rho_p$. In this case, the scaling limit is only performed for the follower population using the arguments  in Section \ref{sub: diffusion limit}. We again obtain \eqref{parfollowers2} for the followers, coupled to the prescribed dynamics of the leaders.  





\section{Inhomogeneous alignment kernel}\label{sec: inhom}
{While \eqref{eq: mean direction} and \eqref{eq: mean direction alt form} are common model equations to relate the mean direction and density, they do not adequately capture key aspects of the follower-leader interaction. For example, because the mean direction $\Lambda$ does not depend on the size, but only the direction of $\mathcal{J}$, a single leader at time $0$ can determine the direction of the swarm.} 

In this section we {therefore also} consider an alignment kernel given by
\begin{equation}
    \Lambda^*=\nu\mathcal{J}(\mathbf{x},t)\ ,\label{eq: new alignemnt}
\end{equation}
where we have removed the normalization and $\mathcal{J}(\mathbf{x},t)$ is given as in (\ref{eq: mean direction}). Here we limit our discussion to the case $n=2$. Parameter $\nu$ is the relaxation frequency, previously assumed constant but in this approach taken to depend on the norm of $\mathcal{J}(\mathbf{x},t)$. Here we study the diffusion and hyperbolic limits of the system for follower-leader interactions, as previously done in Section \ref{sec: macroscopic PDE}. In particular, we focus on the equation for followers (\ref{eq: diffusion followers}) and (\ref{eq: followers hyperb}) under an alignment 
given by (\ref{eq: new alignemnt}). 

The distribution of aligned directions $\Phi(\Lambda^*\cdot\theta)$ in (\ref{eq: alignment}) is replaced by
\begin{equation}
    \bar{\Phi}(\Lambda^*\cdot\theta)=\frac{\Phi(\Lambda^*\cdot\theta)}{\int_0^{2\pi}\Phi(|\Lambda^*|\cos\theta)d\theta}
\end{equation}
such that $\int_S\bar{\Phi}(\Lambda^*\cdot\theta)d\theta=1$.

In the diffusion limit, we replace $\Phi^0(\Lambda\cdot\theta)$ by $\bar{\Phi}^0(\Lambda^*\cdot\theta)$ in equation  (\ref{eq: mean direction followers}). Noting that
\begin{equation}
    \int_S\theta\bar{\Phi}^0(\Lambda^{*}_W\cdot\theta)d\theta=\bar{z}\Lambda^{*}_W\ ,\ \ \textnormal{where}\ \ \ \bar{z}=|\Lambda^*|\int_0^{2\pi}\bar{\Phi}^0(|\Lambda^*|\cos(s))\cos(s) ds\ ,\label{eq: new z}
\end{equation}
for $\theta=\cos(s)\Lambda^*+\sin(s)\Lambda^{*\perp}$ we again obtain (\ref{eq: followers diffusion final}), where in this case we replace $z\Lambda^0_W$ by $\bar{z}\Lambda^*_W$.

For the hyperbolic limit, let us first define the right hand side of (\ref{eq: right hand side}) as
$$
L(\sigma_f)=-\beta\sigma_f+\beta\zeta T\sigma_f+(1-\zeta)\beta\bar{\Phi}^\varepsilon(\Lambda^*\cdot\theta)\rho_f\ .
$$
We know that as $\varepsilon\rightarrow 0$ the solution $\sigma_f^0=\bar{\Psi}(\theta)\rho_f$, where $\bar{\Psi}(\theta)=\zeta+(1-\zeta)\bar{\Phi}^0(\Lambda^*\cdot\theta)$. The new operator $\bar{\Psi}(\theta)$ needs to be a Generalized Collisional Invariant of the operator $L(\sigma_f)$, as in the following sense \cite{degond2013macroscopic,degond2008continuum,dimarco2016self}.
\begin{definition}
A function $\bar{\Psi}(\theta)$ is a Generalized Collisional Invariant of $Q$ if it satisfies
\[
\int_S L(\sigma_f)\bar{\Psi}(\theta)d\theta=0\ ,
\]
for any $\sigma_f$. Equivalently, $\sigma_f$ satisfies 
$
P_\perp\Bigl(\int_S\sigma_f(\mathbf{x},t,\theta)\theta d\theta \Bigr)=0
$, where $P_\perp=\textnormal{Id}-\Lambda^*\otimes\Lambda^*$ is an orthogonal projection to $\Lambda^*$.\label{def: GCI}
\end{definition}

Note that for the case $\zeta=0$, i.e. when only alignment is considered, from Definition \ref{def: GCI} we conclude that $\int_SL(\sigma_f)\bar{\Phi}^0(\Lambda^*\cdot\theta)d\theta=0$, where $\bar{\Phi}^0$ can be taken as the von Mises-Fisher distribution. Then, the analysis in \cite{degond2013macroscopic} follows.
The system (\ref{eq: hyperbolic followers 1})-(\ref{eq: hyperbolic followers 2}) can be written now as
\begin{align}
    \partial_t\rho_f+c_f\bar{z}(1-\zeta)\nabla\cdot(\rho_f\Lambda^*) & =0\ ,\\ \rho_f(\bar{z}(1-\zeta)\partial_t\Lambda^*+\bar{C}_1\Lambda^*\cdot\nabla\Lambda^*)+\bar{C}_2P_\perp\nabla\rho_f & =0\ ,
\end{align}
where $\bar{z}$ is given by (\ref{eq: new z}) and, similar to previous derivation,  $\bar{C}_1=c_f(1-\zeta)\bar{a}_3$, $\bar{C}_2=c_f(1-\zeta)\mathbbm{1}\bar{a}_1+c_f\mathbbm{1}\pi\zeta$ for, $\bar{a}_3=\bar{a}_0-\bar{a}_1$ with
\[
\bar{a}_0=|\Lambda^*|^2\int_0^{2\pi}\bar{\Phi}^0(|\Lambda^*|\cos(s))\cos^2(s)ds\ ,\ \ \bar{a}_1=|\Lambda^*|^2\int_0^{2\pi}\bar{\Phi}^0(|\Lambda^*|\cos(s))\sin^2(s)ds\ .
\]

\section{Discussion}

The capacity of individuals to coordinate their movement to generate collective movement is a phenomenon that has attracted significant interest, in both cellular and animal systems. 
Much of the progress in this area has been facilitated through modelling studies, particularly via the employment of agent-based (or particle) descriptions that consider the movement
response of each single member according to its neighbours. Yet the largest swarms can 
extend over kilometres and contain millions (e.g. herrings, \cite{makris2009}) or even hundreds of 
billions (e.g. desert locusts, \cite{rainey1967,skaf1990}) of members. At such numbers and scales, 
continuous modelling approaches become necessary for their efficiency and increased tractability. Consequently, there is a clear interest in clarifying the relevant form of continuous models 
for swarming systems. 

Much recent interest has focussed on follower-leader systems, where the population is decomposed into a leader population which somehow guides a population of followers, e.g. honey bee swarms where only scouts know the swarm's destination. Inspired by this and similar examples we have formulated a minimal microscopic 
description for a population of followers and leaders, 
deriving the ensuing macroscopic models under distinct scaling limits. To test
the extent to which velocity alignment by itself can propel a guided and coherent 
swarm, follower orientation is limited to the interaction 
choice (\ref{eq: mean direction alt form}): alignment according to
the velocity direction. Under both hyperbolic and parabolic scaling, macroscopic models feature drift-terms with advection in the active leader direction. Thus, alignment to velocity alone yields translocation of the swarm towards the target. 

Under the hyperbolic scaling a pure-drift equation is generated, implying the potential for travelling-pulses with movement of a cohesive colony towards the 
target site. The parabolic limit, on the other hand, generates a drift-diffusion equation. 
While drift is in the direction of the target site, the additional dispersion leads to 
swarm spreading with time. Early stretching of the colony is to be expected, as the initially tight cluster morphs into a migrating swarm. Once established, though, swarms generally retain a relatively stable speed and shape: continued dispersal would be far from optimal, perhaps resulting in population loss or increasing the risk of predation. This, however, must be viewed in light of our intentionally simple modelling approach, where 
we have specifically tested the practicalities of an ``alignment-only'' mechanism. Agent-based approaches for modelling swarms are typically augmented by additional attractive/cohesion behaviours, where individuals are also pulled in the direction of those in their neighbourhood; in non-local continuous models, a similar effect is gained through a nonlocal attraction term that 
biases movement direction, e.g. \cite{bernardi2021}. As stressed above, we have reasonably excluded such considerations 
from the present model for simplicity, however including nonlocal attractions could clearly counteract swarm dispersal. For example, following \cite{bellomo2008modeling}, attraction (or 
repulsion) between individuals can be considered by
\[
\mathcal{A}=\Bigl(\frac{d_c-\textnormal{dist} (\mathbf{x},\mathbf{y})}{\textnormal{dist} (\mathbf{x},\mathbf{y})} \Bigr)e^{-(d_c-\textnormal{dist}(\mathbf{x},\mathbf{y}))^2}\vec{e}(\mathbf{x},\mathbf{y})\ .
\]
Here $\textnormal{dist}(\mathbf{x},\mathbf{y})$ denotes the distance between two individuals $\mathbf{x}$ and $\mathbf{y}$, along the direction of the unit vector  $\vec{e}(\mathbf{x},\mathbf{y})$. If $\textnormal{dist}(\mathbf{x},\mathbf{y})<d_c$, where $d_c$ is a critical distance, then the action of $\mathcal{A}$ is repulsive and if $\textnormal{dist}(\mathbf{x},\mathbf{y})>d_c$ we have attraction.

Recently, there has been considerable interest in the composition and structuring of swarming populations. For example, in the case of bird flocks and fish shoals, the existence of faster and/or braver individuals can lead to hierarchical swarm arrangements 
\cite{nagy2010, reebs2000} and the question is raised as to how much swarm movement is dominated by the choice of a few. For simplicity, the movements of our leaders has been set here
somewhat naively: during streaks, they operate as ballistic particles adopting fast movements towards the target. While this may be a reasonable approximation
for bee swarms, where leaders have specific {\em a priori} knowledge, more general ``leaders''
may be more subtle and get influenced by neighbour movements. Consequently, a 
logical extension would be to also adopt a velocity jump model for the leader population, where one
of their movement contributions stems from an interaction function similar to 
(\ref{eq: mean direction alt form}), distinctly weighted. In other instances, leader/follower
statuses may be transient, for example resulting from spatial position within the swarm, and it
may be necessary to include switching terms between follower and leader populations.

The study of population dynamics, where a few discrete agents act as leaders,
provides an excellent scenario to derive specific control laws and interactions that drive
self-organisation, leading the swarm optimally to a desired outcome.  This
will allow an analytic understanding of follower-leader interactions, with applications
not only to understanding the collective dynamics of biological populations but also outside biology, 
for example to hierarchical swarms of robots. 


\appendix
\section{Macroscopic diffusion limit for followers}\label{sec: Macroscopic diffusion limit}

To compute the mean direction of the followers, $w_f$ in Section \ref{sub: diffusion limit}, we start by substituting the following expansions 
\begin{equation}
    \begin{aligned}
    \sigma_f&=\sigma_f^0+\varepsilon\sigma_f^1+\mathcal{O}(\varepsilon^2)\ ,&& \rho_f=\rho_f^0+\varepsilon\rho_f^1+\mathcal{O}(\varepsilon^2)\ ,\\
    T^\varepsilon& =T^0+\varepsilon T^1+\mathcal{O}(\varepsilon^2)\ ,&& \Phi^\varepsilon(\Lambda\cdot \theta)=\Phi^0+\varepsilon\Phi^1+\mathcal{O}(\varepsilon^2)\ ,\label{eq: expansions followers}
    \end{aligned}
\end{equation}
into \eqref{eq: diffusion followers} and we obtain
\begin{align}
    \varepsilon^2\partial_t(\sigma_f^0+\varepsilon\sigma_f^1)+\varepsilon c_f\theta\cdot\nabla(\sigma_f^0+\varepsilon\sigma_f^1)=&-\beta(\sigma_f^0+\varepsilon\sigma^1_f)+\zeta\beta(T^0+\varepsilon T^1)(\sigma_f^0+\varepsilon\sigma^1_f)\nonumber\\ &+\varepsilon(1-\zeta)\beta(\Phi^0+\varepsilon\Phi^1)(\rho_f^0+\varepsilon\rho_f^1)\ .
\end{align}
Grouping in terms of powers of $\varepsilon$ and assuming that $T^0\sigma_f^0=\rho_f^0$ we have
\begin{align}
   \varepsilon^{0}:\ \ \sigma_f^0 & =\zeta \rho_f^0\ ,\\\varepsilon^{1}:\ \ \sigma_f^1 & =\frac{1}{1-\zeta}\Bigl(-\frac{c_f\zeta}{\beta}\theta\cdot\nabla\rho_f^0+\zeta^2T^1\rho_f^0+(1-\zeta)\Phi^0\rho_f^0 \Bigr)\  .
\end{align}
The Chapman-Enskog expansion in this case is given by
$
\sigma_f=\zeta\rho_f^0+\varepsilon\sigma_f^1+\mathcal{O}(\varepsilon^2)\ .
$
The conservation equation \eqref{eq: conservation equation followers} is given by
\[
\varepsilon^2\partial_t\rho_f^0+ {\varepsilon^2c_f}\nabla\cdot\int_S\theta\sigma_f^1= 0\ .
\]
Finally, we compute the mean direction $w_f$ as follows
\begin{align}
    \nabla\cdot\int_S\theta\sigma_f^1d\theta=-\Delta(D_f\rho_f^0)+\nabla\cdot\rho_f^0\int_S\theta\Phi^0(\Lambda_W\cdot \theta)d\theta\ ,\label{eq: mean direction followers}
\end{align}
where $D_f=\frac{\zeta c_f}{\beta(1-\zeta)}\int_S\theta\theta^Td\theta$. In the derivation of \eqref{eq: mean direction followers} we used $\int_S\theta T^1d\theta=0$ since we assume that $T$ is symmetric.

Also, $W$ is the total mean direction of the whole population and 
\[
\Lambda_W=\frac{\mathcal{J}_W}{|\mathcal{J}_W|}
\ \ \textnormal{for}\ \ 
\mathcal{J}_W={n\varepsilon}\int_\mathbf{y}K^\varepsilon\Bigl(\frac{|\mathbf{y}-\mathbf{x}|}{\varepsilon} \Bigr)W(\mathbf{y},t)d\mathbf{y}\ .
\]
If we consider $\mathcal{J}$ as in (\ref{eq: mean direction alt form}) then the mean direction, $\Lambda$, will depend only on the mean direction of the active leaders, $w_a$, as follows
\[ 
\Lambda_{a}=\frac{\mathcal{J}_a}{|\mathcal{J}_a|}\ \textnormal{where}\ \mathcal{J}_a={n\varepsilon\lambda}\int_{\mathbf{y}}K^\varepsilon\Bigl( \frac{|\mathbf{y}-\mathbf{x}|}{\varepsilon}\Bigr)w_ad\mathbf{y}\ .
\]

The integral over $S$ in (\ref{eq: mean direction followers}) is given by $\int_S\theta\Phi^0(\Lambda^0_W\cdot\theta)d\theta=z\Lambda^0_W$ where $z$ can be computed using polar coordinates $\theta=(\cos(s),\sin(s))$ for $n=2$, or spherical coordinates $\theta=(\cos\phi\sin(s),\sin\phi\sin(s),\cos(s))$ for $n=3$. $z$ is given by \cite{dimarco2016self}
\begin{equation}
  z=\begin{cases}
    \int_{0}^{2\pi}\Phi^0(\cos(s))\cos(s) ds, & \text{if $n=2$},\\
    2\pi\int_{0}^{\pi}\Phi^0(\cos(s))\cos(s)\sin(s) ds, & \text{if $n=3$}.\label{eq: constant z}
  \end{cases}
\end{equation}

\section{Diffusion limit for passive leaders with noise}\label{sec: Diffusion limit for passive leaders}

Let us suppose that the total population of passive leaders, $\sigma_p$, can be written in terms of the following expansion:
\begin{equation}
\sigma_p=\sigma_p^0+\varepsilon\sigma_p^1+\mathcal{O}(\varepsilon^2)\label{eq: passive expansion}\ .
\end{equation}
{Similarly, we expand $\rho_p=\rho_p^0+\varepsilon\rho_p^1+\mathcal{O}(\varepsilon^2)$, $\rho_a=\rho_a^0+\varepsilon\rho_a^1+\mathcal{O}(\varepsilon^2)$} and 
{\begin{equation}\begin{aligned}
B^\varepsilon(\theta)&=B^0(\theta)+\varepsilon B^1(\theta)+\mathcal{O}(\varepsilon^2)\ ,\\
\bar{R}^\varepsilon_{ap}(\theta)&=\bar{R}^0_{ap}(\theta)+\varepsilon\bar{R}^1_{ap}(\theta)+\mathcal{O}(\varepsilon^2)\ ,\ \ \ \ 
{R}^\varepsilon_{ap}={R}^0_{ap}+\varepsilon{R}^1_{ap}+\mathcal{O}(\varepsilon^2)\ ,\label{eq: other expansions}
\end{aligned}\end{equation}}
{where $B^0(\theta)=pT^0(\theta)+(1-p)T^0_{\rho_f}$ and $B^1(\theta)=pT^1(\theta)+(1-p)T^1_{\rho_f}$.}

Substituting  \eqref{eq: passive expansion} and \eqref{eq: other expansions} into (\ref{eq: passive scaled}) and  grouping the appropriate powers of $\varepsilon$ we get
\begin{align}
   \varepsilon^{0}:\ \ \sigma_p^0 & =B^0 \rho_p^0\ ,\label{eq: epsilon-1}\\\varepsilon^{1}:\ \ \sigma_p^1 & =\frac{1}{\beta}(-c_p\theta\cdot\nabla\sigma_p^0{+\beta B^1\rho_p^0+\beta B^0\rho_p^1})\  .\label{eq: epsilon0}
\end{align}
We obtain a  Chapman-Enskog expansion from (\ref{eq: passive expansion}) given by
\begin{equation}
    \sigma_p=B^0\rho^0_p+\varepsilon\sigma_p^1+\mathcal{O}(\varepsilon^2)\ .
\end{equation}
Substituting the above expression into the mean direction of the passive leaders, $w_p$ in (\ref{eq: macroscopic}), we obtain
\begin{align*}
    w_p & =\frac{1}{n}\int_S\theta B^0(\theta)d\theta \rho_p^0+\frac{\varepsilon}{n}\int_S\theta\sigma_p^1d\theta =\frac{\varepsilon}{n}\int_S\theta\sigma_p^1d\theta\ .
\end{align*}
{Since we assume $B^0$ is symmetric, then $\int_S\theta B^0d\theta=0$}. The conservation equation for the passive leaders \eqref{eq: passive scaled2}, is then
\begin{equation}
    \varepsilon^2\partial_t\rho_p^0+{\varepsilon^2  c_p}\nabla\cdot \int_S\theta\sigma_p^1d\theta=\varepsilon^2R_{ap}^0\rho_a^0-\varepsilon^2R_{pa}^0\rho_p^0 \ .
    \end{equation}
The term $\nabla\cdot\int_S\theta\sigma^1_pd\theta$ can be explicitly computed using (\ref{eq: epsilon-1}) and (\ref{eq: epsilon0}) as follows:
\begin{align}
    \nabla\cdot\int_S\theta\sigma_p^1d\theta & =\frac{-c_p}{\beta}\int_S(\nabla\cdot\theta)(\theta\cdot\nabla)B^0\rho_p^0d\theta\nonumber\nonumber\\ & \ \ \ +\frac{1}{\beta}\nabla\cdot\int_S\theta\Bigl({\beta B^1\rho_p^0+\beta B^0\rho_p^1}\Bigr)d\theta\ ,\nonumber\\& =-\Delta (D_p{+\mathbb{T}})\rho_p^0 \ ,\label{eq: mean direction 2}
\end{align}
where we have used the fact that $B(\theta)$ is symmetric and therefore {$\int_S\theta B^0d\theta=\int_S\theta B^1d\theta=0$} and   $$D_p=\frac{pc_p}{\beta}\int_S\theta\theta^TT^0(\theta) d\theta\ ,\ \ \textnormal{and}\ \  {\mathbb{T}=\frac{(1-p)c_p}{\beta}\int_S\theta \theta^T T^0_{\rho_f}d\theta}\  .$$ 


\section{Turn angle operator}\label{sec: turn_angle_properties}
 
This section recalls some basic spectral properties of the turn angle operator $T$ defined in (\ref{eq: turning angle simple}). 
\begin{lemma}\label{lem: eigenfunctions}
Assume that $\tilde{k}$ is continuous. Then $T$ is a symmetric compact operator. In particular, there exists an orthonormal basis of $L^2(S)$ consisting of eigenfunctions of $T$.\\ With $\mathbf{\theta}=(\theta_0,\theta_1,...,\theta_{n-1}) \in S$, we have
\begin{equation}
\begin{aligned}\phi_{0}(\theta) & ={1} &  & \text{is an eigenfunction to the eigenvalue} &  & \nu_{0}=1,\\
\phi_{1}^j(\theta) & ={n\theta_j} &  & \text{are eigenfunctions to the eigenvalue} &  & \nu_{1}=\int_{S}\tilde{k}(\cdot,|\eta-e_1|)\eta_{1}d\eta<1. \label{eq: eigen}
\end{aligned}
\end{equation}
Any function ${\sigma}\in L^2(\mathbb{R}^n\times \mathbb{R}^+\times S)$ admits a unique decomposition 
\begin{equation}
{\sigma}(\mathbf{x},t,\theta)={u}+n\mathbf{\theta}\cdot {w} +\hat{z},\label{eq: real_eigen}
\end{equation}
where $\hat{z}$ is orthogonal to all linear polynomials in $\theta$. Explicitly,
\begin{equation*}
{u}(\mathbf{x},t)=\int_S{\sigma}(\mathbf{x},t,\mathbf{\theta})\phi_0(\theta) d\theta,\ {w}^j(\mathbf{x},t)=\int_S {\sigma}(\mathbf{x},t,\mathbf{\theta})\phi_1^j(\theta) d\theta,
\end{equation*}
and  ${w} =  ({w}^1, \dots, {w}^n)$.
\end{lemma}

\section*{Acknowledgments} 
This research was partially supported by the Italian Ministry of Education, University and Research (MIUR) through the ``Dipartimenti
di Eccellenza'' Programme (2018-2022) -- Dipartimento di Scienze Matematiche ``G. L. Lagrange'', Politecnico di Torino
(CUP: E11G18000350001). SB is member of GNFM (Gruppo Nazionale per la Fisica Matematica) of INdAM (Istituto
Nazionale di Alta Matematica), Italy. KJP acknowledges departmental funding through the “MIUR–Dipartimento di Eccellenza” programme.


\bibliographystyle{abbrv}
\bibliography{bees.bib}

\medskip
Received xxxx 20xx; revised xxxx 20xx.
\medskip

\end{document}